\documentclass[12pt]{article}
\usepackage{amsfonts,amsmath}
\def \beq {\begin{eqnarray}}
\def \eeq {\end{eqnarray}}
\def \beqn {\begin{eqnarray*}}
\def \eeqn {\end{eqnarray*}}

\newcommand{\halmos}{\rule{1ex}{1.4ex}}

\newcounter{for}[section]

\newtheorem{itlemma}{Lemma}[section]
\newtheorem{itproposition}[itlemma]{Proposition}
\newtheorem{theorem}[itlemma]{Theorem}
\newtheorem{itcorollary}[itlemma]{Corollary}
\newtheorem{itremark}[itlemma]{Remark}
\newtheorem{itremarks}[itlemma]{Remarks}
\newtheorem{itdefinition}[itlemma]{Definition}
\newtheorem{itexample}[itlemma]{Example}

\newenvironment{fact}{\begin{itfact}\rm}{\end{itfact}}
\newenvironment{claim}{\begin{itclaim}\rm}{\end{itclaim}}
\newenvironment{lemma}{\begin{itlemma}}{\end{itlemma}}
\newenvironment{remark}{\begin{itremark}}{\end{itremark}}
\newenvironment{remarks}{\begin{itremarks}\rm}{\end{itremarks}}
\newenvironment{corollary}{\begin{itcorollary}}{\end{itcorollary}}
\newenvironment{proposition}{\begin{itproposition}}{\end{itproposition}}
\newenvironment{definition}{\begin{itdefinition}\rm}{\end{itdefinition}}
\newenvironment{example}{\begin{itexample}\rm}{\end{itexample}}
\newenvironment{proof}{\noindent {\em Proof}.\ \
}{\hspace*{\fill}$\halmos$\medskip}
\newcommand{\be}[1]{\addtocounter{for}{1} \begin{equation}\label{#1}}
\newcommand{\ee}{\end{equation}}
\newcommand{\bl}[1]{\begin{lemma}\label{#1}}
\newcommand{\br}[1]{\begin{remark}\label{#1}}
\newcommand{\brs}[1]{\begin{remarks}\label{#1}}
\newcommand{\bt}[1]{\begin{theorem}\label{#1}}
\newcommand{\bd}[1]{\begin{definition}\label{#1}}
\newcommand{\bp}[1]{\begin{proposition}\label{#1}}
\newcommand{\bc}[1]{\begin{corollary}\label{#1}}
\newcommand{\bfact}[1]{\begin{fact}\label{#1}}
\newcommand{\bex}[1]{\begin{example}\label{#1}}
\newcommand{\ec}{\end{corollary}}
\newcommand{\efact}{\end{fact}}
\newcommand{\eex}{\end{example}}
\newcommand{\el}{\end{lemma}}
\newcommand{\er}{\end{remark}}
\newcommand{\ers}{\end{remarks}}
\newcommand{\et}{\end{theorem}}
\newcommand{\ed}{\end{definition}}
\newcommand{\ep}{\end{proposition}}
\newcommand{\epr}{\end{proof}}
\newcommand{\bpr}{\begin{proof}}
\newcommand{\bcl}[1]{\begin{claim}\label{#1}}
\newcommand{\ecl}{\end{claim}}
\newcommand{\ecs}{\end{corollary}}
\newcommand{\eers}{\end{exercise}}
\newcommand{\eexs}{\end{example}}
\newcommand{\eems}{\end{example}}
\newcommand{\els}{\end{lemma}}
\newcommand{\eles}{\end{lemmaex}}
\newcommand{\ets}{\end{theorem}}
\newcommand{\eds}{\end{definition}}
\newcommand{\eps}{\end{proposition}}

\newcommand{\bi}{\begin{itemize}}
\newcommand{\ei}{\end{itemize}}
\newcommand{\ben}{\begin{enumerate}}
\newcommand{\een}{\end{enumerate}}

\def\vbar{\mathchoice{\vrule height6.3ptdepth-.5ptwidth.8pt\kern-.8pt}
   {\vrule height6.3ptdepth-.5ptwidth.8pt\kern-.8pt}
   {\vrule height4.1ptdepth-.35ptwidth.6pt\kern-.6pt}
   {\vrule height3.1ptdepth-.25ptwidth.5pt\kern-.5pt}}
\def\fudge{\mathchoice{}{}{\mkern.5mu}{\mkern.8mu}}
\def\bbc#1#2{{\rm \mkern#2mu\vbar\mkern-#2mu#1}}
\def\bbb#1{{\rm I\mkern-3.5mu #1}}
\def\bba#1#2{{\rm #1\mkern-#2mu\fudge #1}}
\def\bb#1{{\count4=`#1 \advance\count4by-64 \ifcase\count4\or\bba A{11.5}\or
   \bbb B\or\bbc C{5}\or\bbb D\or\bbb E\or\bbb F \or\bbc G{5}\or\bbb H\or
   \bbb I\or\bbc J{3}\or\bbb K\or\bbb L \or\bbb M\or\bbb N\or\bbc O{5} \or
   \bbb P\or\bbc Q{5}\or\bbb R\or\bbc S{4.2}\or\bba T{10.5}\or\bbc U{5}\or
   \bba V{12}\or\bba W{16.5}\or\bba X{11}\or\bba Y{11.7}\or\bba Z{7.5}\fi}}

\def \qed {{\hspace*{\fill}$\halmos$\medskip}}

\def \v {\varphi}
\def \T {{\cal{T}}}

\def \A {{\cal{A}}}
\def \H {{\cal{H}}}

\def \C {{\cal{C}}}
\def \D {{\cal{D}}}

\def \S {{\cal{S}}}

\newcommand{\ba}[1]{\addtocounter{for}{1} \begin{eqnarray}\label{#1}}
\newcommand{\ea}{\end{eqnarray}}

\def\sqr#1#2{{\vcenter{\vbox{\hrule height .#2pt
                             \hbox{\vrule width .#2pt height#1pt \kern#1pt
                                   \vrule width .#2pt}
                             \hrule height .#2pt}}}}

\def\pmb#1{\setbox0=\hbox{#1}%
   \kern-.025em\copy0\kern-\wd0
   \kern.05em\copy0\kern-\wd0
 \kern-.025em\raise.0433em\box0 }
\def\sqr#1#2{{\vcenter{\vbox{\hrule height.#2pt
     \hbox{\vrule width.#2pt height#1pt \kern#1pt
   \vrule width.#2pt}\hrule height.#2pt}}}}

\def\B{{\mathcal B}}
\def\E{{\mathbb E}_0}
\def\P{{\mathbb P}}

\def \J{{\mathcal J}}
\def \II{{\mathcal I}}
\def\N{{\mathbb N}}
\def\Z{{\mathbb Z}}
\def\R{{\mathbb R}}
\def\Q{{\mathcal Q}}

\def\EE{{\mathcal E}}

\def\var{\text{var}}

\def\reff#1{(\ref{#1})}

\def \ind {\hbox{ 1\hskip -3pt I}}
\newcommand {\cro}[1] {\left[ {#1} \right]}
\newcommand {\acc}[1] {\left\{ {#1} \right\}}
\newcommand {\pare}[1] {\left( {#1} \right)}

\newcommand {\bra}[1] {\left< {#1} \right>}

\begin{document}

\title{Annealed upper tails for the energy of a polymer.}
\author{Amine Asselah \\ Universit\'e Paris-Est\\
Laboratoire d'Analyse et de Math\'ematiques Appliqu\'ees\\
UMR CNRS 8050\\ amine.asselah@univ-paris12.fr}
\date{}
\maketitle
\begin{abstract}
We study the upper tails for the energy of a randomly charged symmetric and
transient random walk. 
We assume that only charges on the same site interact pairwise.
We consider {\it annealed} estimates, that is
when we average over both randomness, in dimension three or more.
We obtain a large deviation principle, and an explicit rate function
for a large class of charge distributions.
\end{abstract}



{\em Running head}: Energy upper tails for a polymer.

\section{Introduction} \label{sec-intro}
We consider the following toy model for a charged polymer in dimension
3 or more.
Time and space are discrete, and two independent sources of randomness
enter into the model.
\begin{itemize}
\item[(i)] A symmetric
random walk, $\{S(n),n\in \N\}$, evolving on the
sites of $\Z^d$ with $d\ge 3$.
When the walk starts at $z\in \Z^d$, its law is denoted ${\mathbb P}_z$.
\item[(ii)] A random field of charges, $\{\eta(n), n\in \N\}$.
The charges are centered i.i.d.\,and their law is denoted by $Q$.
We assume that each charge variable satisfies Cramer's condition, that
is for some $\lambda_0>0$, $E_Q[\exp(\lambda_0 \eta)]<\infty$.
\end{itemize}
For a large integer $n$, our polymer is a linear chain
of $n$ monomers each carrying a random charge, and sitting sequentially on
$\{S(0),\dots,S(n-1)\}$. 
The monomers interact pairwise only when they occupy the same site 
on the lattice, and produce a local energy
\be{intro.1}
H_n(z)=\sum_{0\le i\not= j< n} \eta(i)\eta(j) \ind\acc{S(i)=S(j)=z}\qquad
\forall z\in \Z^d. 
\ee
The energy of the polymer, $H_n$, is the sum of $H_n(z)$ over $z\in\Z^d$.

In this paper, we study the annealed probability, i.e.\,when
averaging over both randomness, that $H_n$ is large.
The study of the event $\{-H_n$ large$\}$, that is the lower tails, is
tackled in a companion paper~\cite{A08-lower}, since the phenomenology
and the techniques are different.

A few words of motivation.
Our toy-model comes from physics, where it is used to model
proteins or DNA folding. For mathematical works on
random polymers, we refer to the review paper \cite{HK}, and
the recent books~\cite{giacomin,denH}.

Our interest actually stems from works of Chen~\cite{chen07}, and
Chen and Khoshnevisan~\cite{chen-khosh}, dealing with central limit
theorems for $H_n$ (in the transient, and in the more delicate
recurrent case). The former paper shows some analogy between $H_n$
and the $l_2$-norm of the local times of the walk, whereas the latter
paper shows similarities between the typical fluctuations
of $H_n$ and of a {\it random walk in random scenery}, to be
defined later. Chen~\cite{chen07}
establishes the following annealed moderate deviation principle.
First, some notations: for two positive sequences $\{a_n,b_n,n\in \N\}$, 
we say that $a_n\ll b_n$, when
$\limsup \frac{\log(a_n)}{\log(b_n)}<1$, we
denote the annealed law with $P$, and we denote with $\eta$ a generic
charge variable. 

\bp{prop-chen}[X.Chen~\cite{chen07}] Assume $d\ge 3$, and
$E[\exp(\lambda \eta^2)]<\infty$, for some $\lambda>0$. 
When $\xi_n$ goes to infinity, with $n^{1/2}\ll \sqrt{n}\xi_n\ll n^{2/3}$ ,
then
\be{intro.4}
\lim_{n\to \infty}\  \frac{1}{\xi_n^{2}} \ 
\log\pare{ P(\pm \frac{H_n}{\sqrt n}\ge \xi_n)}= -\frac{1}{2c_d},
\quad\text{where}\quad c_d=\sum_{n\ge 1} \P_0(S(n)=0).
\ee
\ep
Features of the charge distribution do not enter into the moderate
deviations estimates, but play a role in the large deviation regime, which
starts when $\{|H_n|\ge n^{2/3}\}$. Thus, to present our results, we first 
distinguish tail-behaviours of the $\eta$-variables. 
For $\alpha>0$, we say that $\H_\alpha$ holds, or simply that
$\eta\in \H_\alpha$, when $|\eta|^\alpha$
satisfies Cramer's condition.
Also, in order to write shorter proofs, we assume two
non-essential but handy features:
$\eta$ is symmetric with a unimodal distribution (see \cite{AC04}), and
we consider the simplest aperiodic walk: the walk jumps to
a nearest neighbor site or stays still with equal probability.

Finally, we rewrite the energy into a convenient form. For $z\in \Z^d$,
and $n\in \N$, we call $l_n(z)$ the {\it local times},
and $\check{q}_n(z)$ the {\it local charges}. That is
\[
l_n(z)=\sum_{k=0}^{n-1} \ind\acc{S(k)=z},
\quad\text{and}\quad
\check{q}_n(z)=\sum_{k=0}^{n-1} \eta(k) \ind\acc{S(k)=z}.
\]
Inspired by equation (1) of \cite{dgh} and (1.18) of \cite{chen07}, 
we write $H_n(z)=\check{X}_n(z)+Y_n(z)$ with
\[
\check{X}_n(z)=\check{q}_n^2(z)-l_n(z),
\quad \text{and}\quad 
Y_n(z)=l_n(z)-\sum_{i=0}^{n-1}\!\eta(k)^2\ind\acc{S(k)=z}.
\]
Now, 
\be{intro.10}
Y_n=\sum_{z\in\Z^d} Y_n(z)=\sum_{i=0}^{n-1}\pare{1-\eta^2(i)},
\ee
is a sum of centered independent random variables, and its
large deviation asymptotics are well known (see below Remark~\ref{rem-H}). 
Thus, we focus on $\check{X}_n=\sum_{\Z^d} \check{X}_n(z)$. 
\bt{theo-LDP}
Assume $d\ge 3$, and $n^{2/3}\ll \xi_n\ll n^2$. 
If $\eta\in \H_\alpha$ with $\alpha>1$, then
there is a positive constant $\Q_2$,
\be{main-LDP}
\lim_{n\to\infty} \frac{1}{\sqrt{\xi_n}}\log
P\pare{ \check{X}_n \ge \xi_n}= -\ \Q_2.
\ee
Moreover,
the main contribution to $\{ \check{X}_n\ge \xi_n\}$, comes from 
\[
\D_n^*(A)=\{z: A\sqrt{\xi_n}> l_n(z)> \sqrt{\xi_n}/A\},
\]
for some $A>0$. In other words, we have
\be{intro.6}
\limsup_{A\to\infty}\limsup_{n\to\infty} \frac{1}{\sqrt{\xi_n}}
\log\ P\pare{\sum_{z\not\in \D_n^*(A)} \check{X}_n(z) \ge \xi_n}=-\infty.
\ee
\et
Theorem~\ref{theo-LDP} is based ultimately on a subadditive
argument, and the rate $\Q_2$ is beyond the reach of this method. 
However, there is large family of charge distributions for which we can
explicitly compute $\Q_2$. To formulate a more precise result,
we need additional assumptions and notations.

We call the $\log$-Laplace transform of the charge distribution
$\Gamma(x)=\log E_Q[\exp(x\eta)]$. 
Since the charges satisfy Cramer's condition, their empirical
measure obeys a Large Deviation Principle with rate function $\II$,
the Legendre-transform of $\Gamma$:
\be{def-II}
\II(x)=\sup_{y\in \R}\cro{yx-\Gamma(y)}.
\ee
Finally, we define $\chi_d$ (when $d\ge 3$) such that
\[
\P_0(S_k=0,\text{ for some }k>0)=\exp(-\chi_d).
\]
\bt{theo-gauss} Assume $d\ge 3$, and $n^{2/3}\ll \xi_n\ll n^2$. 
Assume $\Gamma$ is twice differentiable and satisfies
\be{hyp-I}
x\mapsto \Gamma(\sqrt x) \text{ is convex on }\R^+.
\ee
Then,
\be{LDP}
\lim_{n\to\infty} \frac{1}{\sqrt{\xi_n}} \log\pare{P(\check{X}_n\ge \xi_n)}=
-\Gamma^{-1}(\chi_d).
\ee
\et
\bex{example-1}
For any $\alpha\in ]1,2]$, we give in Section~\ref{sec-ex}
examples of charge distributions
in $\H_\alpha$ satisfying \reff{hyp-I}. 
\eex
\br{rem-gauss} Note that
when $\eta$ is a centered gaussian variable
of variance $\sigma$, then $\Gamma^{-1}(\chi_d)=
\sqrt{2\chi_d/\sigma}$. In this case, our proof of
Theorem~\ref{theo-gauss} applies equally well
to {\it random walk} in the following {\it random scenery}:
let $\{\zeta(z),z\in \Z^d\}$ be i.i.d.\,with $\zeta(0)$ distributed as 
$Z^2-\sigma$, where $\sigma$ is a positive real, and $Z$ is a centered
gaussian variable with variance $\sigma$. Then, for $d\ge 3$, and
$n^{2/3}\ll \xi_n\ll n^2$, we have
\be{rwrs-H1}
\lim_{n\to\infty}  \frac{1}{\sqrt{\xi_n}} \log\pare{
P(\bra{l_n,\zeta} \ge \xi_n)}= -\sqrt{\frac{2\chi_d}{\sigma}},
\quad\text{where}\quad
\bra{l_n,\zeta}=\sum_{z\in \Z^d} l_n(z)\zeta(z).
\ee
Note that $\zeta\in \H_1$. A large deviations principle for RWRS
was obtained in \cite{ghk} 
in the case $\H_\alpha$ with $0<\alpha<1$
and $\xi_n\gg n^{(1+\alpha)/2}$, and in \cite{A07} with 
\[
1<\alpha<\frac{d}{2},\quad\text{ and}\quad
n^{1-\frac{1}{\alpha+2}}\ll\xi_n \ll n^{1+ \frac{1}{\alpha}}.
\]
To our knowledge, \reff{rwrs-H1} is the first
result for a special case of the border-line regime $\H_1$.
\er
\br{rem-H} Since $H_n=\check{X}_n+Y_n$, and $Y_n$, given in \reff{intro.10},
is a sum of $n$ independent variables bounded by $1$, it is clear 
that \reff{main-LDP}, and \reff{intro.6} hold for 
$H_n$ instead of $\check{X}_n$.
\er
We wish now to present intuitively two ways of understanding 
Theorem \ref{theo-LDP}.
More precisely, we wish to explain why realizing 
an energy of order $\xi_n$, forces a transient polymer to pile of the
order $\sqrt{\xi_n}(\ll n)$ monomers,
on a finite number of sites (independent
of $n$), where the {\it local charge} is of order $\sqrt{\xi_n}$. 
Indeed a first step in proving Theorem \ref{theo-LDP} is the
following result.
\bp{theo1} Assume $d\ge 3$, and $n^{\frac{2}{3}}\ll \xi_n\ll n^2$. 
Consider $\eta$ of type $\H_\alpha$ with $1\le \alpha\le 2$. There are 
positive constants $c_-,c_+$ such that
\be{intro.5}
e^{-c_-\sqrt{\xi_n}}\le P\pare{ \check{X}_n\ge \xi_n}\le
e^{-c_+\sqrt{\xi_n}}.
\ee
When $\xi_n=\xi n^2$, and $\xi>0$ with $Q(\eta>\sqrt\xi)>0$,
\reff{intro.5} holds.
\ep
\reff{intro.6} tells us that the dominant contribution comes from
sites visited $\sqrt{\xi_n}$-times. In other words, a finite number
(independent of $n$)
of short pieces containing of the order of $\sqrt{\xi_n}$ monomers
pile up. Note also that
\[
\sum_{z\in \D_n^*(A)} \check{X}_n(z)=\sum_{z\in \D_n^*(A)} 
\check{q}^2_n(z)+O(\sqrt{\xi_n}),
\]
and we can write
\[
\check{q}^2_n(z)=l_n(z)^2\pare{ \frac{\check{q}_n(z)}{l_n(z)}}^2.
\]
This suggests that on the {\it piles} (i.e.\,sites of $\D_n^*(A)$), the
average charge is of order unity, and the {\it local charges} perform
large deviations.

Assume that $\eta$ is of type $\H_1$. Note that 
$\check{X}_n$ is the $l^2$-norm of an additive random fields
$n\mapsto \{\check{q}_n(z),\ z\in\Z^d\}$.
This is analogous to the self-intersection local times
\[
\|l_n\|^2_2:=\sum_{z\in \Z^d} l_n^2(z).
\]
Now, our (naive) approach in the study of the excess self-intersection 
local time in~\cite{AC04,AC05,A06} is to slice the $l_2$-norm
over the level sets of the local times. By analogy, we
define here the level sets of the {\it local charges}
$\{\check{q}_n(z),\ z\in \Z^d\}$. For a value $\xi>0$
\[
\EE_n(\xi)=\acc{z\in \Z^d:\ \check{q}_n(z)\sim \xi}.
\]
For simplicity, we choose $\xi_n=n^\beta$, and we 
focus on the contribution of $\EE_n(n^x)$ for $0< x\le \beta/2$,
\be{KK.1}
\begin{split}
\acc{\sum_{z\in \EE_n(n^x)} \check{q}_n^2(z)\ge n^\beta}&\subset
\acc{|\EE_n(n^x)|\ge n^{\beta-2x}}\\
&\subset \acc{ \exists
\Lambda\subset[-n,n]^d,\ \&\ |\Lambda|\le n^{\beta-2x}\text{ such that}\
\check{q}_n(\Lambda)\ge |\Lambda| n^{x}\  },
\end{split}
\ee
where $\check{q}_n(\Lambda)$ is the charge collected in $\Lambda$
by the random walk in a time $n$. Thus, \reff{KK.1} requires an 
estimate for $P(\check{q}_n(\Lambda)\ge t)$. Note that by standard estimates,
if we denote by $|\Lambda|$ the number of sites of $\Lambda$,
\be{KK.2}
E[\check{q}_n(\Lambda)^2]=
\sum_{z\in \Lambda} E_Q[\eta^2] \E[l_{n}(z)]\le
\sum_{z\in \Lambda} \E[l_{\infty}(z)]\le
C |\Lambda|^{2/d}.
\ee
\reff{KK.2} motivates the following simple concentration Lemma.
\bl{lem-KK-main}
Assume dimension $d\ge 3$. For some constant $\kappa_d>0$,
and any finite subset $\Lambda$ of $\Z^d$, 
we have for any $t>0$ and any integer $n$
\be{KK-main}
P\pare{ \check{q}_n(\Lambda)\ge t}\le \exp\pare{-\kappa_d \frac{t}
{|\Lambda|^{1/d}}}.
\ee
\el
Note the fundamental difference with the total time spent in $\Lambda$
(denoted $l_\infty(\Lambda)$):
for some positive constants $\tilde \kappa_d$
\be{origin}
P\pare{ l_\infty(\Lambda)\ge t}\le \exp\pare{-\tilde\kappa_d \frac{t}
{|\Lambda|^{2/d}}}.
\ee
We have established \reff{origin} in Lemma 1.2 of \cite{AC04}. 
Thus, using \reff{KK-main} in \reff{KK.1}, we obtain
\be{KK.3}
P\pare{\sum_{z\in \EE_n(n^x)} \check{q}_n^2(z)\ge n^\beta}\le C_n(x)
\exp\pare{-\kappa_d n^{\zeta(x)}},
\ee
with
\be{KK.4}
\zeta(x)=\beta(1-\frac{1}{d})-(1-\frac{2}{d})x,\quad\text{and}\quad
C_n(x)=(2n+1)^{dn^{\beta-2x}}.
\ee
Looking at $\zeta(x)$, we observe that the high level sets (of $\check{q}_n$)
give the dominant contribution and $\zeta(\frac{\beta}{2})=
\frac{\beta}{2}$ in dimension three or more. Note that \reff{KK.4}
also suggests that $d=2$ is a critical dimension, even though
Lemma~\ref{lem-KK-main} fails in $d=2$. 
If one is to pursue this approach rigorously, 
one has to tackle the contribution of $C_n(x)$. Nonetheless, these
simple heuristics show that inequality 
\reff{KK-main} is essentially responsible
for the upper bound \reff{intro.5}.
It is easy to see that \reff{KK-main} is wrong when
$\eta$ is of type $\H_\alpha$ with $0<\alpha<1$, and a different
phenomenology occurs.

First, an observation of Chen \cite{chen07} is that fixing a realization
of the walk, $\{\check{q}_n(z),\ z\in \Z^d\}$ are $Q$-independent random 
variables, and
\be{crucial-chen}
\{\check{q}_n(z),\ z\in \Z^d\}\stackrel{Q-\text{law}}{=}
\{q_n(z),\ z\in \Z^d\}\quad\text{where}\quad
q_n(z)=\sum_{i=1}^{l_n(z)} \eta_z(i),
\ee
where we denote by $\{\eta_z(i),\ z\in \Z^d,\ i\in \N\}$ i.i.d.\,variables
distributed as $\eta$. Also,
\be{anneal-1}
\{\check{X}_n(z),\ z\in \Z^d\}\stackrel{Q-\text{law}}{=}
\{X_n(z),\ z\in \Z^d\}\quad\text{where}\quad
X_n(z)=q_n^2(z)-l_n(z).
\ee
Now, a convenient way of thinking about $X_n$ is to first fix a realization
of the random walk, and to rewrite \reff{anneal-1} as  
\be{intro.7}
X_n(z)= l_n(z) (\zeta_z(l_n(z))-1),\text{where for any $n$}\quad
\zeta_z(n)= ( \frac{1}{\sqrt{n}}\sum_{i=1}^{n} \eta_z(i))^2.
\ee
Thus, $\check{X}_n$ is equal in $Q$-law to a 
scalar product $\bra{l_n,\zeta-1}$
known as {\it random walk in random scenery} (RWRS). However, in our
case the {\it scenery} is a function of the local times.
To make this latter remark more concrete, we recall that when
$\H_\alpha$ holds with $1<\alpha<2$, 
we have some constants $C,\kappa_0,\kappa_\infty$
(see Section~\ref{sec-prel} for more precise statements), such that if
\be{zeta-heuristics}
\zeta(n)=( \frac{1}{\sqrt n} \sum_{i=1}^n \eta(i))^2,\quad
\text{then}\quad Q\pare{\zeta(n)>t}\le C
\left\{ \begin{array}{ll}\exp(-\kappa_0 t)&\mbox{ when }t\ll n\, , \\
\exp(-\kappa_\infty t^{\frac{\alpha}{2}} n^{1-\frac{\alpha}{2}})
& \mbox{ when } n\ll t\, .
\end{array} \right.
\ee
Thus, according to the size of $t/n$, $\zeta(n)$ is either in $\H_1$ or
is a heavy-tail variable (corresponding to $\H_{\alpha/2}$).

Consider the RWRS $\bra{l_n,\rho}$, where
$\{\rho(z),z\in \Z^d\}$ are 
centered independent variables in $\H_\alpha$ with
$\alpha=1+\epsilon$ with a small $\epsilon>0$.
This corresponds to a field of charges with lighter tails than 
$\{\zeta_z(l_n(z)),\ z:\ l_n(z)>0\}$.
Nonetheless, $\{\bra{l_n,\rho}\ge n^\beta\}$ with $\beta>\frac{2}{3}$,
corresponds to a regime (region II of \cite{AC05})
where a few sites, say in a region $\D$, are visited of order 
$\sqrt{\xi_n}$. The phase-diagram of \cite{AC05} suggests
that $\{\bra{l_n,\zeta(l_n)-1}\ge n^\beta\}$ behaves similarly.
Thus, there is a finite region $\D$, over which the sum of 
$\zeta_z(\sqrt{\xi_n})$ should be
of order $\sqrt{\xi_n}$. According to \reff{zeta-heuristics}, this
should make $\{\zeta_z(l_n(z)),z\in \D\}$
of type $\H_1$. It is easy to check, that if
the walk were to spend less time on its most visited sites $\D$, say a time
of order $n^\gamma$ with $\gamma<\frac{\beta}{2}$, then the
$\{\zeta_z(l_n(z)),z\in \D\}$ would be of type $\H_{\alpha/2}$, and an easy
computation using \reff{zeta-heuristics} shows that
\be{heuristics-key}
Q\pare{ \sum_{z\in \D} \zeta_z(n^\gamma)>n^{\beta-\gamma}}\sim 
\exp\pare{ -\kappa_{\infty}n^{\gamma(1-\frac{\alpha}{2})}
n^{(\beta-\gamma)\frac{\alpha}{2}}}\ll \exp\pare{-n^{\frac{\beta}{2}}}.
\ee
This explains intuitively \reff{intro.6}. Note that if $\alpha=1$ in
\reff{heuristics-key}, then for any $\gamma\le \frac{\beta}{2}$, we would
have
$Q(\zeta(n^\gamma)>n^{\beta-\gamma})\sim \exp(-n^{\frac{\beta}{2}})$.

Assume now that the dominant contribution to
the deviation $\{X_n\ge \xi_n\}$ comes from
the random set $\D_n^*(A)$, whose {\it volume is independent of $n$}.
Our next step is to fix a realization of $\D_n^*(A)$, and
integrate over the charges of the {\it monomers piled up in} $\D_n^*(A)$.

Theorem~\ref{theo-gauss} is made possible by the following result.
\bp{lem-gauss}
Assume $\Gamma$ is twice differentiable and satisfies \reff{hyp-I}.
Then, for any finite subset $\D$, any $\gamma>0$, and any 
positive sequence $\{\lambda(z),z\in \D\}$, we have
\be{rate-5}
\inf_{\kappa\ge 0}\cro{\sum_{z\in \D} \lambda(z)  \II\pare{ \kappa(z)}:\ 
\sum_{z\in \D} \lambda^2(z) \kappa^2(z)\ge \gamma^2}=
\pare{\max_{\D} \lambda }\ 
\II\pare{\frac{\gamma}{\max_{\D}\lambda}}.
\ee
Moreover, for any $\alpha,\beta$ positive,
\be{explicit-6}
\inf_{\lambda>0}\cro{ 
\alpha \lambda+ \lambda\II(\frac{\beta}{\lambda})}
=\beta \Gamma^{-1}(\alpha).
\ee
\ep
Without the assumption of Proposition~\ref{lem-gauss}, 
\reff{intro.6} allows us to borrow a strategy developped
in \cite{A07} to prove a large deviation principle for the
self-intersection local times. 
Indeed, the approach of \cite{A07} relies on the fact that
a finite number of piles are responsible for producing the 
excess energy. 

The paper is organized as follows. In Section~\ref{sec-prel}, we 
recall well known bounds on sums of independent random variables,
prove Lemma~\ref{lem-KK-main}, and recall the large deviations
for the $q$-norm of the local times. In Section~\ref{sec-theo1},
we prove Proposition~\ref{theo1}, whose upper bound
is divided in three cases:
$\eta\in \H_2$, in $\H_\alpha$ for $\alpha\in ]1,2[$, and in $\H_1$.
The lower bound in Proposition~\ref{theo1} is established in 
Section~\ref{lower}.
In Section~\ref{sec-gauss}, we prove the large deviation principle
of Theorem~\ref{theo-gauss}. We first discuss useful features of
the rate functions, and then prove Proposition~\ref{lem-gauss}
 in Section~\ref{sec-gauss}. The upper bound of the LDP is
proved in Section~\ref{sec-ubexp}, and the lower bound follows in
Section~\ref{sec-lbexp}.
In Section~\ref{sec-LDP}, we prove Theorem~\ref{theo-LDP}. It is
based on a subadditive argument, Lemma~\ref{lem-subadd}, which
mimics Lemma 7.1 of \cite{A07}. Finally, an Appendix
collects proofs which have been postponed because of their
analogy with known arguments.

\section{Preliminaries}\label{sec-prel}
\subsection{Sums of Independent variables}
In this section, we collect well known results scattered 
in the literature. Since we are not pursuing sharp asymptotics,
we give bounds good enough for our purpose, and  for the convenience
of the reader we have given a proof of the non-referenced results in
the Appendix.

We are concerned with the tail distribution of the $\zeta$-variable,
given in \reff{zeta-heuristics}, with $\bar\zeta(n)=\zeta(n)-1$.
We recall the following result, which we prove in the Appendix
to ease to reading.
\bl{lem-russian}
There are positive constants $\beta_0,\{C_\alpha,\kappa_\alpha,\ 1\le
\alpha\le 2\}$, (depending on the distribution of $\eta$),
such that the following holds.
\begin{itemize}
\item For type $\H_1$, we have
\be{prel.1}
Q\pare{\zeta(n)>t}\le C_1
\left\{ \begin{array}{ll}\exp(-\kappa_1 t)&\mbox{ when }t<\beta_0 n\, , \\
\exp(-\kappa_1 \sqrt{ \beta_0 t n})  & \mbox{ when } t\ge \beta_0 n\, ,
\end{array} \right.
\ee
\item For type $\H_\alpha$, with $1<\alpha<2$, we have
\be{prel.2}
Q\pare{\zeta(n)>t}\le C_\alpha
\left\{\begin{array}{ll}\exp(-\kappa_\alpha t)
&\mbox{ when }t<\beta_0 n\, , \\
\exp\pare{-\kappa_\alpha t^{\alpha/2}(\beta_0 n)^{1-\frac{\alpha}{2}}}
 & \mbox{ when } t\ge \beta_0 n\, ,
\end{array} \right.
\ee
\item For type $\H_2$, we have for any $t>0$,
\be{prel.3}
Q\pare{\zeta(n)>t}\le C_2 \exp(-\kappa_2 t).
\ee
\end{itemize}
\el
A. Nagaev has considered in \cite{anagaev} a sequence $\{\bar Y_n,n\in \N\}$ 
of independent centered i.i.d satisfying $\H_\alpha$ with $0<\alpha<1$, and
has obtained the following upper bound (see also inequality (2.32) of 
S.Nagaev~\cite{snagaev1}).
\bp{prop-anagaev}[of A.Nagaev.] Assume $E[\bar Y_i]=0$ and
$E[(\bar Y_i)^2]\le 1$. There is a constant $C_Y$, such that 
for any integer $n$ and any positive $t$
\be{prel.4}
P\pare{\bar Y_1+\dots+\bar Y_n\ge t}\le C_Y\pare{
nP\pare{\bar Y_1>\frac{t}{2}}+ \exp\pare{-\frac{t^2}{20 n}} }.
\ee
\ep
\br{rem-Y} 
Note that if $\eta\in \H_\alpha$ for $1<\alpha\le 2$, 
then $\eta^2\in \H_{\frac{\alpha}{2}}$. Thus, for $\bar Y_i=\eta(i)^2-1$,
Proposition~\ref{prop-anagaev} yields
\be{term-Y1}
P\pare{\sum_{i=1}^n (\eta(i)^2-1)\ge \xi_n}\le 
C_Y\pare{n\exp\pare{-c_{\alpha}(\xi_n)^{\alpha/2}}+
\exp\pare{ -\frac{\xi_n^2 }{20n}}}.
\ee
\er

When we take $t\sim n^\beta$ in \reff{prel.4}, we have the following 
asymptotical result due to A.Nagaev \cite{anagaev}
(see also Theorem 2.1 and (2.22) of \cite{snagaev1}).
\bl{lem-petrov}
Assume $\{\bar Y_n,n\in \N\}$ are centered independent variables
in $\H_\alpha$ with $0<\alpha<1$. If $\xi_n\gg
n^{\frac{1}{2-\alpha}}$, then for $n$ large enough, we have
\be{prel.5}
P\pare{\bar Y_1+\dots+\bar Y_n\ge \xi_n}\le 2n \max_{k\le n}
P\pare{\bar Y_k>\xi_n}.
\ee
\el
\br{rem-Y2}
With the notations of Remark~\ref{rem-Y}, 
Lemma~\ref{lem-petrov} yields for $\xi_n\gg n^{\frac{1}{2-\alpha/2a}}$
\be{term-Y}
P\pare{\sum_{i=1}^n (\eta(i)^2-1)\ge \xi_n}\le 2n P(\eta^2-1\ge
\xi_n)\le 2C_\alpha n \exp\pare{-c_{\alpha}(\xi_n)^{\alpha/2}}.
\ee
Note that $\alpha>1$ implies that $\frac{1}{2-\alpha/2}>\frac{2}{3}$.

For completeness, note that when $Y_i$-variables are in $\H_1$, 
we can use a special form of
Lemma 5.1 of \cite{A06} to obtain an analogue of \reff{prel.4}.
\er

\subsection{On a concentration inequality}
We prove Lemma~\ref{lem-KK-main}. We assume that for $\lambda_0>0$,
we have $E_Q[\exp(\lambda_0 \eta)]<\infty$.

Note that when $\lambda<\lambda_0/2$, there is a positive constant $C$
such that $E_Q[\exp(\lambda \eta)]\le 1 +C\lambda^2$. Now, note that
\be{con.2}
q_n(\Lambda)
=\sum_{z\in \Lambda}\sum_{i=1}^{l_n(z)} \eta_z(i) 
\stackrel{\text{law}}{=} \sum_{i=1}^{l_n(\Lambda)} \eta_0(i).
\ee
We use Chebychev's inequality, for $\lambda>0$, and integrate only
over the $\eta$-variables
\be{con.1}
Q\pare{ q_n(\Lambda)>t}\le
e^{-\lambda t} 
\pare{E_Q\cro{ e^{\lambda \eta }}}^{l_n(\Lambda)} \le e^{-\lambda t} 
\exp\pare{C\lambda^2 l_n(\Lambda)}.
\ee
Now, using \reff{origin}, if we choose
\be{con.4}
C\lambda^2\le \frac{\tilde \kappa_d}{2|\Lambda|^{2/d}},\quad
\text{then we have}\quad \E\cro{\exp\pare{C\lambda^2 l_n(\Lambda)}}\le 2. 
\ee
Thus, \reff{KK-main} follows at once.

\section{Proof of Proposition~\ref{theo1}.}\label{sec-theo1}
We have divided the proof of the upper bounds in
Proposition~\ref{theo1} into the three cases 
$\H_2$, $\H_\alpha$ and $\H_1$. The lower bound in
\ref{intro.5} is obtained in Section~\ref{lower}.
\subsection{The case $\H_2$.}\label{sec-H2}
We consider first annealed upper bounds for $\{X_n\ge \xi_n\}$.
We choose a sequence $\{\xi_n\}$ such that $\xi_n\gg n^{2/3}$, and
$\epsilon>0$ such that $(\sqrt \xi_n)^{1-\frac{4}{3}\epsilon}\gg
n^{1/3}$.

When averaging first with respect to the charges,
and then with respect to the random walk, we write
\be{intro.9}
P\pare{ X_n\ge \xi_n}=
\E\cro{Q\pare{\sum_{z\in \Z^d} l_n(z) \bar\zeta_z(l_n(z)) \ge \xi_n}}.
\ee
Thus, we think of $X_n$
as a weighted sum of independent centered variables 
$\bar\zeta_z(l_n(z))$. Thanks to the uniform bound \reff{prel.3}, the
dependence of $\zeta_z$ on the local time $l_n(z)$ vanishes.
Indeed, since the $\{\zeta_z\}$ are independent and satisfy Cramer's
condition, we can use Lemma 5.1 of \cite{A06} and obtain
that for some $c_u>0$, any $0<\delta<1$, any finite subset $\Lambda
\subset \Z^d$, any $x>0$, and any integer-sequence $\{k(z),\ z\in \N\}$
\be{intro.13}
Q\pare{ \sum_{z\in \Lambda} \bar\zeta_z(k(z))\ge x}\le
\exp\pare{ c_u |\Lambda| \delta^{2(1-\delta)} \max_n\pare{ E_Q[\zeta(n)^2]}-
\frac{\lambda_1 \delta}{2} x}.
\ee
In order to use \reff{intro.13}, we first need
to uncouple the product $l_n(z)\times\bar \zeta_z$ in \reff{intro.9}.

We fix a constant $A\ge 1$, $b_0=1/A$ and $b_{i+1}=2b_i$.
We define the following subdivision
of $[1,\sqrt{\xi_n}/A]$. We define $i_0$ and $N$ as follows,
$b_{i_0}\le 1<b_{i_0+1}$,
\be{gauss.1}
b_N\le \frac{\sqrt{\xi_n}}{A}<b_{N+1},\quad
\text{we set for } i=i_0,\dots,N, \quad
\D_i=\acc{z:\ b_i\le l_n(z)< b_{i+1}}.
\ee
Furthermore, for $\delta_0$ to be chosen small later (independent of
$A$), we define
\be{gauss.2}
\delta_i=\delta_0\pare{\frac{A b_i }{\sqrt{\xi_n}}}^{(1-\epsilon)}
,\quad
\text{and}\quad p_i=p
\pare{\frac{A b_i }{\sqrt{\xi_n}}}^{\epsilon},\quad
\text{with $p$ such that }\quad\sum_{i\ge 1} p_i=1.
\ee
It is important to note that $p$ is independent on $A$.

Now, we perform the decomposition of $X_n$ in terms of level sets $\D_i$.
Note that for any $A\ge 1$
\be{gauss.4}
\begin{split}
\E\cro{Q\pare{\sum_{z: l_n(z)\le b_{N}}
l_n(z) \bar \zeta_z(l_n(z))\ge \xi_n}}\le&
\sum_{i=i_0}^{N-1} \E\cro{Q\pare{\sum_{z\in \D_i} l_n(z)
\bar \zeta_z(l_n(z))\ge p_i\ \xi_n}}\\
\le &\sum_{i=i_0}^{N-1} Q\pare{\sum_{z\in \D_i} \frac{ l_n(z)}{b_{i+1}}
\bar \zeta_z(l_n(z))\ge \frac{p_i\ \xi_n}{b_{i+1}}}.
\end{split}
\ee
Fix now a realization of the random walk, and let $z\in \D_i$. We have
(recall that $\zeta_z\ge 0$)
\be{gauss.5}
Q\pare{ \frac{ l_n(z)}{{b_{i+1}}}
 \zeta_z(l_n(z))\ge t}\le Q\pare{ \zeta_z(l_n(z))\ge t}\le
e^{-\lambda_1 t}
\ee
We use now Lemma 5.1 of \cite{A06}, that we have recalled in \reff{intro.13}
with the choice of $\delta_i$ given in \reff{gauss.2} (and whose dependence
on $n$ is omitted)
\be{gauss.6}
\begin{split}
Q&\pare{\sum_{z\in \D_i} \frac{ l_n(z)}{b_{i+1}}
\bar \zeta_z(l_n(z))\ge \frac{p_i\ \xi_n}{b_{i+1}}}\\
&\qquad\qquad\le\exp\pare{ c_0 |\D_i|\delta_i^{2(1-\delta_i)}-
\frac{\delta_i p_i \xi_n}{4 b_{i+1}}},
\quad\text{with}\quad c_0=c_u \sup_k E_Q[\zeta^2(k)].
\end{split}
\ee
If we denote $\kappa_0=\delta_0^{-2\delta_0}$, then note that
$\delta_i^{-2\delta_i(n)}\le \kappa_0$.
Then, the bound \reff{gauss.6} is useful if the first term on the right hand
side is negligible, that is if
\be{gauss.7}
8 \kappa_0 c_0 |\D_i| \delta_i^2\le \delta_i \ p_i\ 
\frac{\xi_n}{b_{i+1}}.
\ee
Note that
\be{gauss.3}
\delta_i p_i\ \frac{\xi_n}{b_{i+1}}\ge p\frac{Ab_i}{\sqrt{\xi_n}}\ 
\frac{\xi_n}{b_{i+1}}=\frac{p}{2}A \sqrt{\xi_n}.
\ee
Assuming \reff{gauss.7}, the result follows right away by \reff{gauss.3}.
The remaining point is to show that \reff{gauss.7} holds. Note that
\reff{gauss.7} holds as long as $|\D_i|$ is not {\it large}. On the other
hand, we express $\{|\D_i|$ {\it large}$\}$ as a large deviation event
for the self-intersection local time. Thus, \reff{gauss.7} holds when
\be{gauss.8}
\begin{split}
8 \kappa_0 c_0 |\D_i|\le&\frac{1}{2\delta_i^2} p\ A \sqrt{\xi_n}=
\pare{\frac{\sqrt{\xi_n}}{Ab_i}}^{2-2\epsilon} \frac{p}{2}
A \sqrt{\xi_n}\\
\le & \frac{(\sqrt{\xi_n})^{3-2\epsilon}}
{b_i^{2-2\epsilon}} \frac{p}{2A^{1-2\epsilon}}.
\end{split}
\ee
When \reff{gauss.8} does not hold, we note for some constant $c_1$,
(independent of $A$) and any $q>1$
\be{gauss.9}
\EE_i= \acc{ |\D_i|> \frac{c_1 }{A^{1-2\epsilon}}
 \frac{(\sqrt{\xi_n})^{3-2\epsilon}}{b_i^{2-2\epsilon}}}\subset
\acc{ \| l_n\|^q_q\ge \frac{c_1}{A^{1-2\epsilon}}
(\sqrt{\xi_n})^{3-2\epsilon} b_i^{q-2+2\epsilon}}.
\ee
Note that the event on the right hand side of \reff{gauss.9} is
a large deviation event since from our hypotheses on $\xi_n$ and
$\epsilon$, we have, when $q\ge 2$ and in the worse case where $b_i=1$, that
\[
\E[\| l_n\|^q_q]\sim \kappa(q,d)
n\ll \frac{c_1}{A^{1-2\epsilon}}(\sqrt{\xi_n})^{3-2\epsilon}.
\]

\noindent{{\underline{Case $d=3$.}}}

We consider \reff{gauss.9} with $d=3$ and $q=2<q_c(3)=3$, and
use Theorem 1.1 and Remark 1.3 of \cite{A08}. This yields
\be{gauss.10}
\begin{split}
\P(\EE_i)\le &
\P\pare{\| l_n\|_2^2-E\cro{\| l_n\|_2^2}\ge 
\frac{c_1}{2A^{1-2\epsilon}}(\sqrt{\xi_n})^{3-2\epsilon}
b_i^{2\epsilon}}\\
\le& C \exp\pare{ -c(2,3) \pare{ 
\frac{c_1}{2A^{1-2\epsilon}}(\sqrt{\xi_n})^{3-2\epsilon}
b_i^{2\epsilon}}^{2/3}n^{-1/3}}.
\end{split}
\ee
Thus, the quantity $\sum_i \P(\EE_i)$ is negligible if
\be{gauss-neglige}
\sqrt{\xi_n}n^{1/3}\ll 
\pare{\frac{c_1}{2A^{1-2\epsilon}}(\sqrt{\xi_n})^{3-2\epsilon} }^{2/3}
\ee
Condition \reff{gauss-neglige} is equivalent to $n^{2/3}\ll \xi_n$.

\noindent{{\underline{Case of $d\ge 4$.}}}
First, in dimension 5 or more, we have from Theorem 1.2 of \cite{A08}, 
with $q=2>q_c(d)$
\be{gauss.11}
\P(\EE_i)\le
\exp\pare{ -c(2,d) \pare{
\frac{c_1}{2A^{1-2\epsilon}}(\sqrt{\xi_n})^{3-2\epsilon}
b_i^{2\epsilon}}^{1/2}}.
\ee
This term is negligible.

When $d=4$, $q_c(4)=2$, and we choose $q>2$. 
\be{gauss.d4}
\P(\EE_i)\le
\exp\pare{ -c(q,d) \pare{
\frac{c_1}{2A^{1-2\epsilon}}(\sqrt{\xi_n})^{3-2\epsilon}
b_i^{q-2+2\epsilon}
}^{1/q}}
\ee
This term is negligible if
\[
\sqrt{\xi_n}\ll(\sqrt{\xi_n})^{3/q} \Leftarrow\quad q<3.
\]
Thus, if we choose $2<q<3$, $P(\EE_i)\ll \exp(-\sqrt{\xi_n})$.

In conclusion, we obtain that
\be{gauss.12}
\begin{split}
P\big(\sum_{l_n(z)<\sqrt{\xi_n}/A}& X_n(z)\ge \xi_n\big)\le
\sum_{i=i_0}^{N-1} \E Q
\pare{ \sum_{z\in \D_i} l_n(z) \bar \zeta_z\pare{l_n(z)}\ge p_i\xi_n}\\
&\qquad\qquad\le  \sum_{i=i_0}^{N-1}
\E Q\pare{ \sum_{z\in \D_i} \bar \zeta_z\pare{l_n(z)}\ge
p_i\frac{\xi_n}{b_{i+1}} ,\ |\D_i|<
\frac{c_1}{A^{1-2\epsilon}}\frac{(\sqrt{\xi_n})^{3-2\epsilon}}
{b_i^{2-2\epsilon}}}
\\
&\qquad\qquad\quad+\sum_{i=i_0}^{N-1}
\P\pare{ |\D_i|\ge 
\frac{c_1}{A^{1-2\epsilon}}\frac{(\sqrt{\xi_n})^{3-2\epsilon}}
{b_i^{2-2\epsilon}}}\\
&\qquad\qquad\le C\log(n)\times
\exp\pare{ -\frac{p}{8}A\sqrt{ \xi_n}}+
\exp\pare{-\xi_n^{1/2+\epsilon}}.
\end{split}
\ee
Thus, taking $A=1$ in \reff{gauss.12}, we cover the levels
$\{z:\ l_n(z)< \sqrt{\xi_n}\}$, whereas taking $A$
larger than 1, we cover the levels $\{z:\ l_n(z)< \sqrt{\xi_n}/A\}$
Thus, combining these two regimes, we obtain the upper bound \reff{intro.5},
whereas taking $A$ to infinity, we obtain the asymptotic \reff{intro.6}.

\subsection{The case $\H_\alpha$ with $1<\alpha<2$.}\label{sec-alpha}
Charges in $\H_\alpha$ have a much fatter tails than in $\H_2$.
Thus, we decompose $\zeta$ into its small
and large values. For $z\in \Z^d$, define for all positive integer $k$
\be{trunc-zeta}
\zeta_z'(k)= \zeta_z(k)\ \ind\acc{\zeta_z(k)\le \beta_0 k}\quad
\text{and}\quad
\zeta_z''(k)= \zeta_z(k)\ \ind\acc{\zeta_z(k)> \beta_0 k}.
\ee
We add a bar on top of $\zeta',\zeta''$ to denote the centered variables,
and we define 
\[
\bar X_n'=\sum_{z\in \Z^d} l_n(z)\bar\zeta_z'(l_n(z)),\quad\text{and}\quad
\bar X_n''=\sum_{z\in \Z^d} l_n(z)\bar\zeta_z''(l_n(z)).
\]
Note that
\be{Xn-decomp}
\acc{\bar X_n\ge \xi_n}\subset \acc{\bar X_n'\ge \frac{\xi_n}{2}}\cup
\acc{\bar X_n''\ge \frac{\xi_n}{2}}.
\ee
The $\{\zeta_z',\ z\in \Z^d\}$ look like coming from $\eta$ in $\H_2$.
Indeed, for any $t>0$
\be{small.1}
\{\zeta_z'(k)>t\}=\{t<\zeta_z(k)\le \beta_0 k\}\Longrightarrow 
Q\pare{ \zeta_z'(k)>t}\le C_\alpha\exp(-\kappa_\alpha t).
\ee
Thus, the term $\{\bar X_n'\ge n^\beta \frac{\xi}{2}\}$ follows the
same treatment as that of Section~\ref{sec-H2}, with the upper bound
\reff{intro.5}, and the asymptotic \reff{intro.6}.

We focus on the large values of $\zeta_z$. Note that
\be{large.1}
\{\zeta_z''(k)>t\}=\{\zeta_z(k)\ge \max(t,\beta_0 k)\}\Longrightarrow
\forall t>0\quad
Q\pare{ \zeta_z''(k)\ge t}\le C_\alpha e^{
-\kappa_\alpha t^{\frac{\alpha}{2}}(\beta_0 k)^{1-\frac{\alpha}{2}}}.
\ee
For convenience, set $\tilde \alpha=\frac{2}{\alpha}-1$, with
$0<\tilde \alpha<1$, and note that \reff{large.1} implies that
for $u>0$
\be{large.2}
Q\pare{k^{\tilde\alpha}  \zeta_z''(k)\ge u}\le C_\alpha \exp\pare{
-\kappa_\alpha \beta_0^{1-\frac{\alpha}{2}} u^{\frac{\alpha}{2}}}.
\ee
We can therefore think of 
\[
Y_z:=(l_n(z))^{\tilde\alpha}  \zeta_z''(l_n(z)),\quad
\pare{\bar Y_z:=(l_n(z))^{\tilde\alpha}\bar\zeta_z''(l_n(z))}
\]
as having a heavy-tail (of type $\H_{\frac{\alpha}{2}}$). Using the level
decomposition of Section~\ref{sec-H2}, we first fix a realization
of the random walk and estimate 
\be{large.3}
\A_i:=Q\pare{\sum_{z\in \D_i}
\pare{\frac{l_n(z)}{b_{i+1}}}^{1-\tilde \alpha}
\bar Y_z\ge p_i\ \frac{\xi_n}{b_{i+1}^{1-\tilde \alpha}}}.
\ee
Note that from \reff{large.2}, we have
some constant $C$ such that for $z\in \D_i$
\be{large.4}
Q\pare{ \pare{\frac{l_n(z)}{b_{i+1}}}^{1-\tilde \alpha}
Y_z\ge u}\le Q(Y_z\ge u)\le C_\alpha \exp(-C u^{\frac{\alpha}{2}}).
\ee
This implies that for some $\sigma_Y>0$, we have 
$E[Y_z^2]\le \sigma_Y^2$, and we can use Proposition~\ref{prop-anagaev}
\be{large.5}
\A_i\le C_Y\pare{|\D_i| 
Q\pare{Y_1\ge \frac{p_i\ \xi_n}{2(b_{i+1})^{1-\tilde \alpha}}}
+\exp\pare{ -\frac{p_i^2}{20\sigma_Y^2 |\D_i|}
\pare{\frac{\xi_n}{b_{i+1}^{1-\tilde \alpha}}}^2}}.
\ee
We show now that the first term of the right hand side of \reff{large.5}
is the dominant term. Note that 
\be{large.6}
Q\pare{Y_1\ge \frac{p_i\ \xi_n}{2(b_{i+1})^{1-\tilde \alpha}}}\le
C_\alpha\exp\pare{ -C 
\pare{\frac{p_i\ \xi_n}{2(b_{i+1})^{1-\tilde \alpha}}}^{\alpha/2}}.
\ee
Now, we estimate the exponent in the right hand side of \reff{large.6}
\be{large.9}
\frac{p_i\ \xi_n}{2(b_{i+1}})^{1-\tilde \alpha}=p
\pare{\frac{Ab_{i+1}}{2\sqrt{\xi_n}}}^{\epsilon} \frac{\xi_n}
{{2(b_{i+1}})^{1-\tilde \alpha}}.
\ee
When $\tilde\alpha<1$, we choose $\epsilon$ small enough
so that $1-\tilde\alpha>\epsilon$, and then 
the least value of the last term in \reff{large.9} is for $i$ such that
$b_{i+1}=\sqrt{\xi_n}/A$. Thus
\be{large.10}
\frac{p_i\ \xi_n}{2(b_{i+1})^{1-\tilde \alpha}}\ge \frac{p}{2^\epsilon}
\frac{2A^{1-\tilde \alpha}}{(\sqrt{\xi_n})^{1-\tilde \alpha}}\xi_n.
\ee
When taking a power $\alpha/2$ in \reff{large.10},
the power of $\xi_n$ in the right hand side of \reff{large.10}
satisfies
\be{large.7}
\frac{\alpha}{2}(1-\frac{1-\tilde \alpha}{2})=\frac{1}{2}
\Longrightarrow
\pare{\frac{p_i\ \xi_n}{2(b_{i+1})^{1-\tilde \alpha}}}^{\alpha/2}\ge
\pare{\frac{p}{2^{1+\epsilon}}}^{\alpha/2} 
A^{\frac{\alpha}{2} (1-\tilde \alpha)} \sqrt{\xi_n}.
\ee
To deal with $|\D_i|$ in the
second term on the right hand side of \reff{large.5}
(the gaussian bound), 
we can assume as in Section~\ref{sec-H2} that we restrict
ourselves to $\EE_i$ defined in \reff{gauss.9}. Thus, on $\EE_i$
\be{large.8}
\frac{p_i^2}{|\D_i|}\pare{\frac{\xi_n}{(b_{i+1})^{1-\tilde \alpha}}}^2\ge
\frac{p_i^2}{c_1}A^{1-2\epsilon}\xi_n^{\frac{1}{2}+\epsilon}
b_{i+1}^{2(\tilde \alpha-\epsilon)}=\frac{p^2 A}{c}
b_{i+1}^{2\tilde \alpha}\sqrt{\xi_n}\ge \frac{p^2 A}{c} \sqrt{\xi_n}.
\ee
\reff{large.8} holds as soon as $n$ is large enough. 
Thus, the gaussian bound is is negligible.
The proof is concluded as in Section~\ref{sec-H2}.

\qed
\subsection{The case $\H_1$ .}\label{sec-H1}
In this case, no information on the dominant level set can be obtained.
We perform the same decomposition as in Section~\ref{sec-alpha}.
The term $X_n'$ follows the case $\H_2$, and we focus on $X_n''$.
Note that $\tilde \alpha=1$, so that the analysis of
Section~\ref{sec-alpha} is not adequate (note that \reff{large.9}
would yield an upper bound of the form $\exp(-\xi_n^{1\epsilon})$). 
In our case,
\be{large.20}
Y_z=l_n(z)\zeta_z''(l_n(z)),\quad\text{and}\quad
\bar Y_z=Y_z-E[Y_z].
\ee
Note that \reff{prel.1} implies that on $\{l_n(z)>0\}$ and for
some constant $c_1$ we have
$t>0$,
\be{large.21}
\begin{split}
Q(Y_z>t)=& Q\pare{ \zeta_z(l_n(z))>
\max\pare{ \beta_0 l_n(z),\frac{t}{l_n(z)}}}\\
\le & C_1 \exp\pare{ - \kappa_1 \sqrt{ \beta_0  l_n(z)
\max\pare{ \beta_0 l_n(z),\frac{t}{l_n(z)}}}}\\
\le & C_1 \exp\pare{ - \kappa_1\max\pare{ \beta_0 l_n(z),
\sqrt{ \beta_0 t}}}\le C_1 e^{-c_1 \sqrt t}.
\end{split}
\ee
When fixing a realization of the random walk, we think of
$Y_z$ as in $\H_{\frac{1}{2}}$, and use Lemma~\ref{lem-petrov}
with $\xi_n\gg n^{\frac{2}{3}}$, and $n$ large to obtain
\be{large.22}
\begin{split}
P(X_n''\ge \xi_n)=&\E\cro{Q\pare{\sum_{z:l_n(z)>0}
\bar Y_z\ge \xi_n}}\le 2\E\cro{n\max_{z:l_n(z)>0}
Q(\bar Y_z\ge \xi_n)}\\
\le & 2nC_1 \exp\pare{- c_1 \sqrt{ \xi_n}}.
\end{split}
\ee
\subsection{Lower bound in \reff{intro.5}}\label{lower}
A scenario compatible with the cost in \reff{intro.5} is as follows.
The walk is pinned at the origin a time $t_n$ of order $\sqrt{\xi_n}$,
building up an energy $t_n \bar \zeta(t_n)$ required to be
of order $\xi_n$. The remaining
time, the walk roams freely, and the total energy should be made of
$t_n \bar \zeta(t_n)$, and a part close to zero due to the 
centering. Note that the relevant order for $\zeta(t_n)$ is $t_n\sim
\sqrt{\xi_n}$, thus, we are in the {\it central limit regime}. 
This is why we can treat at once the three cases we have considered.

We first show a more general lemma.
\bl{lem-LB} Let $\{\D_n,\ n\in \N\}$ be a sequence of random subsets, with
$\D_n\in ]-n,n[^d$ and $\D_n$ measurable with respect to
$\sigma(S(k),k<n)$. Let $\{m_n, M_n,\ n\in \N\}$
be positive sequences with
$m_n\le n$. Assume either (i) $\xi_n=\xi n^2$ and $Q(\eta>\sqrt{\xi})>0$,
or (ii) $n^{2/3}\ll \xi_n\ll n^2$. Then, for any $\epsilon>0$, 
\be{LB.aim-bis}
\begin{split}
P\pare{ \|q_n\|_2^2 -n \ge \xi_n}&\ge
2^{-2M_n-1} P\pare{\|\ind_{\D_n}q_{m_n}\|_2^2\ge (1+\epsilon) \xi_n
,\ |\D_n|\le M_n}\\
&-\P_0\pare{l_n(\D_n)\ge \frac{\epsilon}{2} \xi_n,\ |\D_n|\le M_n}-
\P_0\pare{\|l_n\|_2\ge \frac{\xi_n}{n^{\epsilon}}}.
\end{split}
\ee
\el
Based on Lemma~\ref{lem-LB}, we distinguish two cases: (i) $\xi_n=\xi n^2$
with $\xi<1$, and (ii) $\xi_n \ll n^2$. 

First, we state a corollary of Theorems 1.1 and 1.2 (and Remark 1.3)
of \cite{A08} whose immediate proof is omitted.
\bc{cor-2beta}
Assume $d\ge 3$ and $\xi_n\gg n^{2/3}$. For $\epsilon>0$ small enough,
and $n$ large enough
\be{neglect-2beta}
\P_0\pare{ \|l_n\|_2 \ge \frac{\xi_n}{n^{\epsilon}}}\le
\exp\pare{ -\sqrt{\xi_n} n^{\epsilon}}.
\ee
\ec

In case (i), we choose $\epsilon$ small enough so that
$Q(\eta>\sqrt{(1+\epsilon)\xi})>0$.
Our scenario is obtained as we choose $\D_n=\{0\}$, $M_n=1$ and
$m_n=n$ in Lemma~\ref{lem-LB}. Inequality \reff{LB.aim-bis}
reads 
\be{scena-0}
P\pare{ \|q_n\|_2^2 -n \ge \xi n^2}\ge
\frac{1}{4} Q\pare{ \frac{1}{n}\sum_{i=1}^{n} \eta_0(i)\ge
1}-\P_0\pare{l_n(0)\ge \frac{\epsilon}{2} \xi n^2}-
\P_0\pare{\|l_n\|_2\ge n^{2-\epsilon}}.
\ee
Now, using $\{\eta_0(i)\ge \sqrt{(1+\epsilon)\xi},\ \forall i\le n\}\subset
\{\eta_0(1)+\dots+\eta_0(n)\ge \sqrt{(1+\epsilon)\xi}n \}$, as well as 
\reff{origin}, and Corollary~\ref{cor-2beta}, we have
\be{scena-2}
P\pare{ \|q_n\|_2^2 -n \ge \xi_n}\ge
\frac{1}{8} Q(\eta>\sqrt{(1+\epsilon)\xi})^{n}.
\ee
In case (ii), we take any $\xi>0$, and 
$m_n^2=(1+\epsilon)\xi_n$. Then, \reff{LB.aim-bis} reads
\be{scena-1}
\begin{split}
P\pare{ \|q_n\|_2^2 -n \ge \xi_n}\ge&
\frac{1}{4} Q\pare{ \sum_{i=1}^{m_n} \eta_0(i)\ge
m_n}-\P_0\pare{l_n(0)\ge \frac{\epsilon}{2} \xi_n}-
\P_0\pare{\|l_n\|_2\ge \frac{\xi_n}{n^{\epsilon}}}\\
\ge & \frac{1}{4} Q\pare{ \eta\ge 1}^{m_n}-
\P_0\pare{l_n(0)\ge \frac{\epsilon}{2} \xi_n}-
\P_0\pare{\|l_n\|_2\ge \frac{\xi_n}{n^{\epsilon}}}.
\end{split}
\ee
To deal with the two last terms in \reff{scena-1}, we use
\reff{origin}, and Corollary \reff{cor-2beta}.

\noindent{\bf Proof of Lemma~\ref{lem-LB}.}

To establish \reff{LB.aim-bis}, we need (i)
to control the process outside $\D_n$, and (ii)
to control the process in the time period $[m_n,n[$.
We start with (i), and introduce
notations $\S_n=\{|\D_n|\le M_n\}$, and 
$\B_n= \acc{\|l_n\|_2\le \xi_n n^{-\epsilon}}$.
First,
\be{LB-subset}
\begin{split}
\acc{\|q_n\|^2_2-n\ge \xi_n}\supset&
\acc{\|\ind_{\D_n}\ q_n\|^2_2-l_n(\D_n)
\ge(1+\frac{\epsilon}{2})\xi_n}\\
&\qquad\qquad\cup \acc{\|\ind_{\D_n^c}\ q_n\|^2_2-l_n(\D_n^c)\ge
-\frac{\epsilon}{2}\xi_n}.
\end{split}
\ee
Note that if we show that 
\be{LB.5}
\ind_{\B_n}\ Q\pare{ \|\ind_{\D_n^c} q_n\|^2_2- l_n(\D_n^c)
\le -\frac{\epsilon}{2} \xi_n}\le \frac{1}{2}\ind_{\B_n},
\ee
then, using the independence of the charges on different regions,
we would have that
\be{LB.6}
2 \ind_{\B_n}
Q\pare{\|q_n\|^2_2-n\ge \xi_n,\ \S_n}
\ge \ind_{\B_n}
Q\pare{\|\ind_{\D_n}q_n\|^2_2-l_n(\D_n)\ge(1+\frac{\epsilon}{2})\xi_n,
\ \S_n}.
\ee
Then, upon integrating \reff{LB.6} over the random walk, we would reach
\be{LB.4}
\begin{split}
P\pare{ \|q_n\|_2^2-n \ge \xi_n}\ge &
P\pare{ \B_n\cap \S_n,\ \|q_n\|_2^2-n \ge \xi_n}\\
\ge &\frac{1}{2} P\pare{ \B_n\cap \S_n,\
\|\ind_{\D_n}q_n\|^2_2-l_n(\D_n)\ge (1+\frac{\epsilon}{2})\xi_n}\\
\ge & \frac{1}{2}P\pare{\|\ind_{\D_n}q_n\|^2_2-l_n(\D_n)\ge
(1+\frac{\epsilon}{2})\xi_n,\ \S_n}-\frac{1}{2}P(\B_n^c).
\end{split}
\ee
We now show \reff{LB.5}. We expand $q_n^2$
\be{LB.9}
q_n^2(z)-l_n(z) =\big(\sum_{i\le l_n(z)}\eta_z(i)\big)^2-l_n(z)
=\sum_{i\le l_n(z)}(\eta^2_z(i)-1)
+2\sum_{1\le i<j\le l_n(z)} \eta_z(i)\eta_z(j),
\ee
It is immediate to obtain, for $\chi_1=E[\eta^4]+1$
\be{LB.10}
E_Q\cro{(q_n^2(z)-l_n(z))^2}
= l_n(z)\pare{E_Q[\eta^4]-1}+2\pare{l_n^2(z)-l_n(z)}\le
\chi_1 l_n^2(z).
\ee
By Markov's inequality 
\be{LB.11}
\begin{split}
\ind_{\B_n}\ Q\pare{\|\ind_{\D_n^c} q_n\|^2_2- l_n(\D_n^c)\le -\epsilon 
\xi_n}\le
&\ind_{\B_n}\ \frac{\sum_{z\not\in \D_n} \var(q_n^2(z)-l_n(z))}
{(\epsilon \xi_n)^2}\\
&\le \ind_{\B_n}\ 
\frac{\chi_1 \|\ind_{\D_n^c}l_n\|^2_2}{(\epsilon \xi_n)^2}
\le \ind_{\B_n}\ \frac{\chi_1}{\epsilon^2}n^{-\epsilon}.
\end{split}
\ee
Thus, for any $\epsilon>0$, \reff{LB.5} holds for $n$ large enough.

We now deal with (ii), and show that
\be{LB.1}
\begin{split}
P\big(\|\ind_{\D_n} q_n\|_2^2-l_n(\D_n)\ge& 
(1+\frac{\epsilon}{2}) \xi_n, \S_n\big)
+\P_0\pare{l_n(\D_n)\ge \frac{\epsilon}{2} \xi_n,\ \S_n}\\
&\ge \pare{\frac{1}{2}}^{2M_n}
P\pare{\|\ind_{\D_n} q_{m_n}\|_2^2\ge (1+\epsilon) \xi_n,\ \S_n}.
\end{split}
\ee
We impose that the {\it local charge} on
each $z\in \D_n$ during both time periods $[0,m_n[$ and
$[m_n,n[$ be of a same sign. Indeed, this would have the effect that
\be{LB.2}
\forall z\in \D_n\quad
q^2_{[0,m_n[}(z)+q^2_{[m_n,n[}(z)
\le \pare{ q_{[0,m_n[}(z)+q_{[m_n,n[}(z)}^2.
\ee
Thus, if we set
\be{LB.3}
\tilde \S_n= \acc{\|\ind_{\D_n} q_{m_n}\|_2^2\ge (1+\epsilon) \xi_n}
\cap \acc{\forall z\in \D_n,\ q_{[0,m_n[}(z)\ge 0},
\ee
and,
\be{LB-time}
\S_n'=\acc{\forall z\in \D_n,\ q_{[m_n,n[}(z)\ge 0},
\quad\text{then}\quad
\tilde \S_n\cap \S_n'
\subset \acc{\|\ind_{\D_n} q_n\|^2_2\ge (1+\epsilon) \xi_n}.
\ee
Note also that when integrating over the charges,
$\tilde \S_n$ and $\S_n'$ are independent, and by symmetry of the
charges' distribution
\be{LB-sym}
Q\pare{\tilde \S_n}\ge \pare{\frac{1}{2}}^{|\D_n|}
Q\pare{\S_n},
\quad\text{and}\quad
Q\pare{\S_n'}\ge  \pare{\frac{1}{2}}^{|\D_n|}.
\ee
Now, \reff{LB-time} and \reff{LB-sym} imply that
\be{LB-(i)}
P\pare{\|\ind_{\D_n} q_n\|^2_2\ge (1+\epsilon)\xi_n ,\ \S_n}\ge
\pare{\frac{1}{2}}^{2M_n}
P\pare{\|\ind_{\D_n} q_{m_n}\|^2_2\ge (1+\epsilon) \xi_n,\ \S_n}.
\ee
Now we center $\|\ind_{\D_n} q_n\|^2_2$. Note that
\be{LB-center}
\acc{\|\ind_{\D_n} q_n\|^2_2\ge (1+\epsilon) \xi_n}\subset
\acc{\|\ind_{\D_n} q_n\|^2_2- l_n(\D_n)\ge(1+\frac{\epsilon}{2}) \xi_n}
\cup \acc{ l_n(\D_n)\ge \frac{\epsilon}{2} \xi_n}.
\ee
Now \reff{LB.1} follows from \reff{LB-(i)} and \reff{LB-center}.

In conclusion, \reff{LB.aim-bis} is obtained as we put together
\reff{LB.4}, and \reff{LB.1}.

\qed

\section{Explicit rate functions}\label{sec-explicit}
\subsection{Some examples}\label{sec-ex}
For constants $0<a, \beta$, and normalizing constant $c(a,\beta)$,
consider the {\it even} density (on $\R$)
\be{ex-3}
g(x)=c(a,\beta) \int_0^1  \exp\pare{-\frac{a}{u^\beta}-u x^2} du.
\ee
When $x$ is large, the value $u^*$ where $\frac{1}{u^\beta}+u x^2$
reaches a minimum,
arises as you equal $\frac{a\beta}{u^\beta}=u x^2$. Thus,
\be{ex-4}
\lim_{x\to\infty} \frac{1}{x^\alpha}
\log g(x)=-c_\alpha, \quad
\text{with}\quad \alpha=\frac{2\beta}{1+\beta}, \quad
\text{and}\quad
c_\alpha=\frac{(1+\beta)}{\beta^{\beta/(1+\beta)}}
a^{1/(1+\beta)}.
\ee

Now, the log-Laplace transform is
\be{ex-5}
\begin{split}
\Gamma(y)=& \log c(a,\beta)\int_0^1 du \int_{-\infty}^\infty dx
\pare{\exp\pare{yx-\frac{a}{u^\beta} -u x^2}}\\
=&\log\pare{c(a,\beta)
 \int_0^1du\ \exp\pare{\frac{y^2}{4u}-\frac{a}{u^\beta}}
\times \int_{\R} \exp\pare{-\pare{x\sqrt u -y/(2\sqrt u )}^2} dx}\\
=&\log \int_0^1du\ \exp\pare{\frac{y^2}{4u}}\phi(u),\quad
\text{with}\quad \phi(u)=c(a,\beta)
\exp(-\frac{a}{u^\beta}) \frac{1}{\sqrt u}\\
\end{split}
\ee
The last integral converges if $\beta>1$, that is $\alpha\in ]1,2]$.
Note, that at infinity $\Gamma(y)\sim d_\alpha y^{2\beta/(\beta-1)}$, as
expected by Tauberian theorems. Now,
$\Gamma(\sqrt y)$ is convex, using H\"older's inequality on the
representation of the last time of \reff{ex-5}.

\subsection{When $x\mapsto \II(\sqrt{x})$ is concave.}\label{sec-gauss}
In this section, we review some useful property of the rate
function, and prove Proposition~\ref{lem-gauss}.
First, we state a simple observation.
\bl{lem-Gamma-I} Assume that $\Gamma$ is twice differentiable, and
$y\mapsto \Gamma(\sqrt y)$ is convex for $y>0$.
Then, $x\mapsto \II(\sqrt x)$ is concave for $x>0$. 
\el
\bpr
Note first that $\Gamma$ is strictly convex. Indeed, note
first that $\Gamma(x)>0$
for $x\not= 0$. This forces $\Gamma'(x)>0$ for $x>0$.
Second, $y\mapsto \Gamma(\sqrt y)$ convex, for $y>0$, implies that
$\Gamma'(x)/x$ is increasing for $x>0$, which in turn says that
$\Gamma'(x)$ is strictly increasing, which implies that $\Gamma$ is
strictly convex.

Also, $\II$ is differentiable, and $\II'$ is the
inverse of $\Gamma'$ on $\R^+$. Note that $\Gamma''(x)>0$ for
$x>0$, so that $\II'$ is differentiable and
$\II''(x) \Gamma''(\II'(x))=1$. 

Now, $\Gamma'(x)/x$ is increasing for $x>0$ is equivalent to 
\[
\forall y>0,\quad\Gamma'(y)\le y \Gamma''(y)\Longrightarrow
\forall x>0,\quad\II'(x)\ge x \II''(x).
\]
Thus, $x\mapsto \II(\sqrt x)$ is concave for $x>0$.
\epr

\noindent{\bf Proof of Proposition~\ref{lem-gauss}.}
We show first two useful properties.
First, note that for $0\le p\le 1$, and $x> 0$
\be{I-convex}
p \II(x)\ge \II(px),
\ee
with equality if and only if $p=1$.
Indeed, using that $\II(0)=0$, \reff{I-convex} is equivalent to
\be{rate-increase}
\frac{\II(x)-\II(0)}{x}\ge \frac{\II(px)-I(0)}{px},\quad\text{with}\quad
0\le xp\le x.
\ee
The strict convexity of $\II$ implies that \reff{rate-increase} is true.
Secondly, for any $x_1,\dots,x_n$ positive
\be{I-concave}
\II(\sqrt{x_1})+\dots+\II(\sqrt{x_n})\ge \II\pare{ \sqrt{x_1+\dots+x_n}}.
\ee
It is easy to see that \reff{I-concave} is obtained by induction as a direct
consequence of $x\mapsto \II(\sqrt x)$ concave and $\II(0)=0$.

Now, assume that $\sum_{z\in \D} \lambda^2(z) \kappa^2(z)
\ge \gamma^2$, and let
$z^*\in \D$ be such that $\lambda(z^*) =\max_{z\in \D} \lambda(z)$.
By using that $\II$ is increasing in $\R^+$ (\reff{I-convex} is stronger
than this latter property)
\be{rate-1}
\sum_{z\in \D} \frac{\lambda^2(z) \kappa^2(z)}{\lambda^2(z^*)}\ge 
\frac{ \gamma^2}{\lambda^2(z^*)}\Longrightarrow
\II\pare{ \sqrt{ \sum_{z\in \D} \frac{\lambda^2(z) \kappa^2(z)}
{\lambda^2(z^*)}}}
\ge \II\pare{ \frac{\gamma}{\lambda(z^*)}}.
\ee
By \reff{I-concave}, we have
\be{rate-2}
\sum_{z\in \D} \II\pare{ \frac{\lambda(z) }{\lambda(z^*)}\kappa(z)}\ge
\II\pare{ \frac{\gamma}{\lambda(z^*)}}.
\ee
From \reff{I-convex}, we deduce that
\be{rate-4}
\sum_{z\in \D} \lambda(z)  \II\pare{ \kappa(z)}\ge
\lambda(z^*)\ \II\pare{ \frac{\gamma}{\lambda(z^*)}}.
\ee

Note that the inequality in \reff{rate-4} is an equality if and only if
$\kappa(z)=0$ for all $z\not= z^*$, and $\lambda(z^*)\kappa(z^*)=\gamma$. 
Thus, \reff{rate-5} holds.

We prove now \reff{explicit-6}.
First, note that since $\II$ is differentiable
\be{explicit-2}
\inf_{x>0}\cro{ \alpha x+ x\II(\frac{\beta}{x})}=
\alpha x^*+ x^*\II(\frac{\beta}{x^*}),
\ee
where $x^*$ satisfies
\be{explicit-3}
\alpha=-\II(\frac{\beta}{x^*})+\frac{\beta}{x^*}\II'(\frac{\beta}{x^*})
\quad(\text{and}\quad \alpha x^*+ x^*\II(\frac{\beta}{x^*})=
\beta\II'(\frac{\beta}{x^*})).
\ee
Recall now that for any $x$
\be{explicit-4}
-\II(x)+x\II'(x)=\Gamma(\II'(x)).
\ee
Thus, combining \reff{explicit-3} and \reff{explicit-4}, we obtain
\be{explicit-5}
\alpha=\Gamma\pare{\II'(\frac{\beta}{x^*})}.
\ee
The parenthesis in \reff{explicit-3}, and \reff{explicit-5} imply 
\reff{explicit-6}, and $x^*$ satisfies
\be{explicit-7}
\frac{\beta}{x^*}=\Gamma'\pare{\Gamma^{-1}(\alpha)}.
\ee
\qed
\subsection{A corollary of Proposition~\ref{theo1}.}
We first derive a corollary of Proposition~\ref{theo1}.
\bc{lem-D}
For any $\epsilon>0$, there are positive constants
$A$ and $M_A$ such that, for $n$ large enough
\be{D-1}
P\pare{ X_n\ge \xi_n}\le 2
P\pare{\|\ind_{\D_n^*(A)} q_n\|_2^2\ge (1-\epsilon)\xi_n
,\ |\D_n^*(A)|\le M_A},
\ee
and,
\be{D-2}
\frac{1}{4^{M_A+1}}P\pare{\|\ind_{\D_n^*(A)} q_n\|_2^2\ge 
(1-\epsilon)\xi_n,\ |\D_n^*(A)|\le M_A}\le P\pare{ X_n\ge \xi_n}.
\ee
\ec
Since, $\epsilon$ is eventually taken to 0,
Theorem~\ref{theo-gauss} follows from Corollary~\ref{lem-D}, once we
find the same upper and lower bound for
\be{ldp.1}
P\pare{
\|\ind_{\D_n^*(A)} q_n\|_2^2\ge \xi_n,\ |\D_n^*(A)|\le M_A}.
\ee

\noindent{\bf Proof of Corollary~\ref{lem-D}.}

Fix $\epsilon>0$ small enough, and $A,M_A$ to be chosen later. 
\be{gauss-5}
\begin{split}
P(X_n\ge \xi_n)\le & P\pare{ \sum_{z\in \D_n^*(A)} q_n^2(z)\ge
(1-\epsilon)\xi_n, |\D_n^*(A)|\le M_A}+R_n\\
\text{with}\qquad & R_n=P(|\D_n^*(A)|\ge M_A)+
P\pare{\sum_{z\not\in \D_n^*(A)} X_n(z)\ge \epsilon \xi_n}.
\end{split}
\ee
From \reff{intro.6}, there is $A$, such that
\be{gauss-2}
P\pare{\sum_{z\not\in \D_n^*(A)} X_n(z)\ge \epsilon \xi_n}
\le \frac{1}{4} e^{-c_-\sqrt{\xi_n}}.
\ee
Also, an application of \reff{origin} stated as Lemma 2.2. of \cite{AC05}
shows that
\be{gauss-3}
P(|\D_n^*(A)|\ge M_A)\le
P\pare{ |\{z:\ l_n(z)\ge \frac{\sqrt{\xi_n}}{A}\}|\ge M_A}
\le |B(n)|^{M_A} \exp(-\tilde \kappa_d \frac{M_A^{1-2/d}\sqrt{\xi_n}}{A}).
\ee
Thus, there is $M_A$ such that
\be{gauss-4}
P(|\D_n^*(A)|\ge M_A)\le
 \frac{1}{4} e^{-c_-\sqrt{ \xi_n}}.
\ee
Now, for $A$ large enough, and the corresponding $M_A$ such that
\reff{gauss-4} holds, we have
\be{gauss-6}
R_n\le \frac{1}{2} \exp\pare{-c_-\sqrt{ \xi_n}}.
\ee
Now, from the lower bound in \reff{intro.6}, we have \reff{D-1}. 

We turn now to \reff{D-2}. We invoke
Lemma~\ref{lem-LB} with $\D_n=\D_n^*(A)$ and $M_n=M_A$, and $m_n=n$.
Note that from \reff{origin}, we have
\be{D-3}
\P_0\pare{l_n(B(r))\ge \frac{\epsilon}{2} \xi_n,\ 
|\D_n^*(A)|<M_A}\le (2n)^{dM_A}
\exp\pare{-\frac{\tilde \kappa_d \epsilon \xi_n}{2M_A^{2/d}}},
\ee
which is negligible, as well as the term $P(\|l_n\|_2\ge
\xi_n n^{-\epsilon})$ by Corollary~\ref{cor-2beta}.
\qed
\subsection{Upper Bound in Proposition~\ref{lem-gauss}.}\label{sec-ubexp}
Our first task is to approximate $\|\ind_{\D_n^*(A)}q_n\|_2$
by a convenient discrete object. 
\paragraph{Step 1: On discretizing the local charge.}
First, write for any integer $n$
\be{gauss-7}
\|\ind_{\D_n^*(A)}q_n\|^2_2=\sum_{z\in \D_n^*(A)} 
l_n^2(z)\pare{\frac{q_n(z)}{l_n(z)}}^2.
\ee
Note now that for any $\epsilon>0$,
there is $\delta>0$ such that on $\{|\D_n^*(A)|\le M_A\}$,
\be{gauss-8}
\acc{\|\ind_{\D_n^*(A)}q_n\|_2^2\ge \xi_n}\subset
\acc{\sum_{z\in \D_n^*(A)} l_n^2(z)
\pi_{\delta}\cro{\pare{\frac{q_n(z)}{l_n(z)}}^2}
\ge (1-\epsilon)\xi_n}.
\ee
Indeed, 
\be{gauss-9}
\pi_{\delta}\cro{\pare{\frac{q_n(z)}{l_n(z)}}^2}\ge
\pare{\frac{q_n(z)}{l_n(z)}}^2-\delta.
\ee
We sum \reff{gauss-9} over $z\in \D_n^*(A)$, on the event $
\{|\D_n^*(A)|\le M_A\}$, and choose $\delta$ small enough so
that $\delta A^2M_A\le \epsilon$, and
\be{gauss-10}
\begin{split}
\sum_{z\in \D_n^*(A)} l_n^2(z)\pi_{\delta}\cro{
\pare{\frac{q_n(z)}{l_n(z)}}^2}\ge& \|\ind_{\D_n^*(A)}q_n\|^2_2-\delta
\sum_{\D_n^*(A)} l_n^2(z)\\
\ge&  \|\ind_{\D_n^*(A)}q_n\|^2_2-\delta A^2M_A \xi_n\\
\ge& (1-\delta A^2M_A)\xi_n\ge (1-\epsilon)\xi_n.
\end{split}
\ee
\paragraph{Step 2: Integrating over the charges.}
As usual, we integrate first with respect to the $\eta$-variables.
We introduce the following random set
(of volume independent of $n$ since $|\D_n^*(A)|\le M_A$)
\[
\B_n=\acc{{\bf \kappa}
=\{\kappa(z),z\in \D_n^*\}: \kappa^2(z) \in \delta \N,\ 
0\le \kappa(z)\le 2A^2\xi_n}.
\]
Now, for $z\in \D_n^*(A)$, $q_n(z)/l_n(z)$ satisfies a Large
Deviation Principle with rate function $\II$:
for $z\in \D_n^*(A)$ (recalling that this implies that
$l_n(z)\ge \sqrt{\xi_n}/A$), $\epsilon(n)$ vanishing as $n$
goes to infinity, and $\kappa(z)>0$
\be{ldp-eta}
\begin{split}
Q\pare{ \pi_\delta\cro{(\frac{q_n(z)}{l_n(z)})^2}=\kappa^2(z)}&\le
2Q\pare{\frac{1}{l_n(z)}\sum_{i=1}^{l_n(z)}\eta_i\ge \kappa(z)}\\
&\le 2\exp\pare{ -l_n(z)\pare{ \II\pare{\kappa(z)}+\epsilon(n)}}.
\end{split}
\ee

Also, if we denote 
\be{def-Cn}
\C_n(l_n)=\acc{{\bf \kappa}
\in \B_n:\ \sum_{z\in \D_n^*(A)} l_n^2(z)\kappa^2(z)\ge 
(1-\epsilon)\xi_n},
\ee
we have, when integrating only over the charge variables,
on the event $\{|\D_n^*(A)|\le M_A\}$,
\be{ldp.4}
\begin{split}
Q\big( \sum_{z\in \D_n^*(A)} & l_n^2(z) 
\pi_\delta\cro{(\frac{q_n(z)}{l_n(z)})^2}\ge (1-\epsilon)\xi_n\big)\\
&\le \pare{ \frac{A\sqrt{\xi_n}}{\delta}}^{|D_n^*(A)|}
\sup_{{\bf \kappa}\in \C_n(l_n)} \prod_{z\in \D_n^*(A)} Q\pare{  
\pi_\delta\cro{(\frac{q_n(z)}{l_n(z)})^2}= \kappa^2(z)}\\
&\le \pare{ \frac{A\sqrt{\xi_n}}{\delta}}^{M_A}
\sup_{{\bf \kappa}\in \C_n(l_n)} \exp\pare{-\sum_{z\in \D_n^*(A)}
l_n(z)\pare{ \II\pare{\kappa(z)}+\epsilon(n)}}\\
&\le \pare{ \frac{A\sqrt{\xi_n}}{\delta}}^{M_A} e^{\epsilon(n)\sqrt{\xi_n}}
\exp\pare{-\inf_{{\bf \kappa}\in \C_n(l_n)}
\sum_{z\in \D_n^*(A)} l_n(z) \II\pare{\kappa(z)}}.
\end{split}
\ee
\paragraph{Step 3: On an explicit infimum.}
We apply Proposition~\ref{lem-gauss} to the infimum in \reff{ldp.4},
since we take actually the infimum over a smaller (discrete) 
set $\C_n(l_n)$.
\be{annealed-1}
\ind_{\{|\D_n^*(A)|\le M_A\}}
Q\big( \sum_{z\in \D_n^*(A)}  q_n^2(z)\ge\xi_n \big)
\le e^{\epsilon(n) \sqrt{\xi_n}}
\exp\pare{ -\max_{\D_n^*(A)} l_n\ \times \II\pare{
\frac{\sqrt{(1-\epsilon)\xi_n}}{\max_{\D_n^*(A)} l_n}}}.
\ee
We now integrate over the random walk (over the event 
$\{|\D_n^*(A)|\le M_A\}$)
\be{ldp.6}
\begin{split}
E\big(\|\ind_{\D_n^*(A)}&q_n\|_2^2\ge\xi_n,\ |\D_n^*(A)|\le M_A\big)\\
&\le e^{\epsilon(n) \sqrt{\xi_n}}
\E\cro{\exp\pare{ -\max_{\D_n^*(A)} l_n \times 
\II\pare{\frac{\sqrt{(1-\epsilon)\xi_n}}{\max_{\D_n^*(A)} l_n}}}}\\
&\le e^{\epsilon(n) \sqrt{\xi_n}} \sum_{k\ge 1} 
\P_0\pare{\max_{z\in\Z^d} l_n(z)=k} 
\exp\pare{ -k \II\pare{\frac{\sqrt{(1-\epsilon)\xi_n}}{k}}}\\
&\le e^{\epsilon(n) 
\sqrt{\xi_n}} \sum_{z\in \Z^d} \P_0(H_z\le n) \sum_{k\ge 1} 
\P_0\pare{ l_\infty(z)=k}
\exp\pare{ -k \II\pare{\frac{\sqrt{(1-\epsilon)\xi_n}}{k}}}.
\end{split}
\ee
Now, a simple coupling argument shows that for any $z$, $\P_0(l_\infty(z)=k)
\le \P_0(l_\infty(0)=k)=\exp(-\chi_d k)$. 
Now, for $\epsilon>0$ arbitrarily small, we call $x=k/\sqrt{\xi_n}$, and
we use \reff{explicit-6}, of Proposition~\ref{lem-gauss}, 
with $\alpha=(1-\epsilon)\chi_d$, and $\beta=\sqrt{(1-\epsilon)}$
for $\epsilon$ small. It is clear from \reff{explicit-7}
that we can choose $A$ large enough so that $x^*\in [1/A,A]$. Thus,
\be{explicit-1}
\inf_{1/A\le x\le A}\cro{(1-\epsilon)\chi_d x+
x\II\pare{\frac{\sqrt{(1-\epsilon)}}{x}}}=
\sqrt{(1-\epsilon)}\times \Gamma^{-1}\pare{(1-\epsilon)\chi_d}.
\ee
Since the function we optimize is continuous, it is irrelevant
whether we take $x$ real in $[1/A,A]$, or along a subdivision of
mesh $1/\sqrt{\xi_n}$ as $n$ goes to infinity. 
\subsection{Lower Bound in Proposition~\ref{lem-gauss}.}\label{sec-lbexp}
Recall Corollary~\ref{lem-D}, and \reff{D-2}. Note that
\[
P\big(\sum_{z\in \D_n^*(A)} q_n^2(z)\ge \xi_n\big)
\ge P\pare{q_n^2(0)\ge \xi_n,\ \{0\}\in \D_n^*(A)}.
\]
When $A$ is large enough (recall that $\beta<2$)
\be{ldp.8}
P\pare{q_n^2(0)\ge\xi_n,\ \{0\}\in \D_n^*(A)}
\ge \sup_{A\sqrt{\xi_n}\ge m_n \ge \sqrt{\xi_n}/A}\P_0\pare{ l_n(0)=m_n}
2Q(\sum_{i=1}^{m_n}\eta_0(i)\ge \sqrt{ \xi_n}).
\ee
We first need to compare $\P_0\{ l_n(0)=m_n\}$ with 
$\P_0\{ l_\infty(0)=m_n\}$, where $m_n=\lfloor x\sqrt{\xi_n}\rfloor$,
(the integer part of $x\sqrt{\xi_n}$), and $x\in [1/A,A]$.
We state the following lemma, which we prove at the end of the section.
\bl{lem-time}
Assume $d\ge 3$. For any $\epsilon>0$, there is $t(\epsilon)>0$, such
that for $n\ge t(\epsilon) m_n$ (with any $m_n\le A\sqrt{\xi_n}$)
\be{ldp.9}
\P_0\pare{ l_\infty(0)=m_n}\le e^{\epsilon m_n} \P_0\pare{ l_n(0)=m_n}.
\ee
\el
Now, recalling that for any integer $k$,
$\P_0(l_\infty(0)=k)=\exp(-\chi_d k)$, $\beta<2$,
and using Lemma~\ref{lem-time}, we have
(with $\epsilon(n)$ is a vanishing sequence, and the supremum
over $m_n$ in $[\sqrt{\xi_n}/A, A \sqrt{\xi_n}]$)
\be{ldp.10}
\begin{split}
P\pare{q_n^2(0)\ge\xi_n,\ \{0\}\in \D_n^*(A)}
\ge & e^{-\epsilon m_n}
\sup_{m_n} \exp\pare{-\chi_d m_n
-m_n\II\pare{\frac{\sqrt{\xi_n}}{m_n}}-\epsilon(n)m_n}\\
\ge & e^{-(\epsilon+\epsilon(n))m_n}
\exp\pare{-\sqrt{\xi_n}\inf_{A\ge x\ge 1/A}
\pare{\chi_d x +x\II\pare{\frac{1}{x}}}}\\
\ge & e^{-(\epsilon+\epsilon(n))m_n}
\exp\pare{-\sqrt{\xi_n}\times \Gamma^{-1}\pare{\chi_d}} .
\end{split}
\ee
We used in the second line of \reff{ldp.10} the continuity of 
the infimum.
Also, we need to choose $A$ large enough so that $x^*$ which
minimizes the infimum in \reff{ldp.10} is in $[1/A,A]$.
A lower bound identical to the upper bound follows 
from \reff{ldp.10}, as we send $\epsilon$ to zero.

\noindent{\bf Proof of Lemma~\ref{lem-time}.}
Recall that $m_n=\lfloor x\sqrt{\xi_n}\rfloor$, and $x\in [1/A,A]$.
Let $\{\tau_i,i\ge 1\}$ be the successive return times to 0, and
recall the classical bound, which holds in $d\ge 3$ for some constant $c_d$
\be{ldp.11}
\P(\tau_i>t|\tau_i<\infty)\le \frac{c_d}{t^{d/2-1}}.
\ee
\reff{ldp.11} implies that for any $\epsilon$, there is $t(\epsilon)$
such that
\be{ldp.12}
\prod_{i\le m_n} \P(\tau_i\le t(\epsilon)|\tau_i<\infty)\ge\pare{
1- \frac{c_d}{t(\epsilon)^{d/2-1}}}^{m_n}\ge \exp(-\epsilon m_n).
\ee
Also, note that 
\be{ldp.13}
\acc{l_\infty(0)=m_n}=\A_n\cap\acc{\tau_{m_n+1}=\infty}
\quad\text{with}\quad
\A_n=\acc{\tau_i<\infty,\ \forall i=1,\dots,m_n}.
\ee
Now
\be{ldp.14}
P(\A_n)=P\pare{ \sum_{i=1}^{m_n} \tau_i<m_n t(\epsilon)|\A_n}
P(\A_n)+P\pare{ \sum_{i=1}^{m_n} \tau_i\ge m_n t(\epsilon)|\A_n}
P(\A_n).
\ee
We show that the first term on the right hand side of \reff{ldp.14}
is large enough. Using \reff{ldp.12}
\be{ldp.15}
P\pare{ \sum_{i=1}^{m_n} \tau_i<m_n t(\epsilon)|\A_n}\ge
\prod_{i\le m_n} P\pare{ \tau_i< t(\epsilon)|\tau_i<\infty}\ge
\exp(-\epsilon m_n).
\ee
Thus, \reff{ldp.14} and \reff{ldp.15} yield
\be{ldp.16}
P(\A_n)\le P\pare{ \sum_{i=1}^{m_n} \tau_i<m_n t(\epsilon)|\A_n}
P(\A_n)+\pare{1-\exp(-\epsilon m_n)}P(\A_n),
\ee
and this implies that
\be{ldp.17}
P(\A_n)P(\tau_{m_n+1}=\infty)\le e^{\epsilon m_n} 
P\pare{ \sum_{i=1}^{m_n} \tau_i<m_n t(\epsilon)}P(\tau_{m_n+1}=\infty).
\ee
Note that \reff{ldp.17} is equivalent to \reff{ldp.9} when $n\ge
m_n t(\epsilon)$.
\qed

\section{A general large deviation principle.}\label{sec-LDP}
To prove Theorem~\ref{theo-LDP}, we
follow the approach of \cite{A07}. We recall
the main steps of the approach, and we detail how to treat the
features which are different. 
In all of Section~\ref{sec-LDP}, we assume that dimension is 3 or more.

In \cite{A07}, dimension is 5 or more. There is actually three
occurrences in \cite{A07} where $d>3$ is used, and we now review them.
\begin{itemize}
\item For self-intersection local times, $d=4$ is the
critical dimension, and only for $d>4$, do we have that
the excess self-intersection is made up on a finite number
of sites. Here, the phenomenology is different, with Lemma~\ref{lem-KK-main}
suggesting that $d=2$ is critical, and
Proposition~\ref{theo1} holds for $d\ge 3$.
\item For normalizing time in Section 6 of \cite{A07}, we used
that, conditionned on returning to 0,
the return time to 0 has finite expectation if $d>4$. We bypass
this constraint in Section~\ref{sec-UBLDP} 
(see the arguments following \reff{time-1}).
\item Lemma 4.9 of \cite{A07} uses an estimate on the probability
of exiting a sphere from a given domain in (10.29). This latter
inequality is useful in $d>3$. In $d=3$, one can use instead the
more sophisticated estimate of Lemma 5 (b) of \cite{LBG} which 
states that for $z\in B(r)$, and $\Sigma$ a domain on the boundary
of $B(r)$, the probability a random walk starting on $z$
exits $B(r)$ in $\Sigma$ is bounded by a constant time $|\Sigma|/
|z-r|^{d-1}$ (rather than $|\Sigma|/|z-r|^{d-2}$). 
Thus, the denominator of (10.30) of \cite{A07}
has a power $2d-3$ (rather than $2d-4$), 
and (10.31) holds also in $d=3$ once
$L$ is chosen large enough, where $L$ is related to the diameter
of a ball containing the piles of monomers producing the
excess energy (see $\tilde \Lambda$ in the paragraph 
following \reff{sub.2infinity}).
\end{itemize}
\subsection{On a subadditive argument}\label{sec-sub}
We recall that Lemma 7.1 of \cite{A07} establishes that for any
radius $r$, and $\xi>0$, there is a positive constant
$\J(\xi,r)$, and the following limit exists
\[
\lim_{n\to\infty} 
\frac{1}{n}\log\ \P_0\pare{\|\ind_{B(r)}l_n\|_2\ge \xi n,\ S(n)=0}
=-\J(\xi,r).
\]
One important difference between local times, and local charges, is
that the distribution of the latter is continuous. Thus, we cannot
find an optimal strategy by
maximizing over a finite number of values, as for $\{
\|\ind_{B(r)}l_n\|_{2}\ge \xi n\}$.
The remedy is to first discretize $\|\ind_{B(r)}q_m\|_{2}$.

For any $\epsilon>0$, there is $\delta>0$ such that
\be{sub.8}
\acc{\|\ind_{B(r)}q_m\|_2\ge \xi m}\subset
\acc{\sum_{z\in B(r)} l_m^2(z)
\pi^2_{\delta}\pare{\frac{q_m(z)}{l_m(z)}}
\ge ((1-\epsilon)m)^2}.
\ee
Indeed, if $z$ is such that $q_m(z)\ge l_m(z)$, then
\be{sub.9bis}
\pi_{\delta}\pare{ \frac{q_m(z)}{l_m(z)}}
\ge \frac{q_m(z)}{l_m(z)}-\delta\ge (1-\delta)\frac{q_m(z)}{l_m(z)}.
\ee
Now, if $q_m(z)<l_m(z)$, then
\be{sub.9}
\pi^2_{\delta}\pare{\frac{q_m(z)}{l_m(z)}}\ge
\pare{\frac{q_m(z)}{l_m(z)}}^2-2\delta.
\ee
We use now \reff{sub.9bis} and \reff{sub.9} to form $\|q_m\|^2_{B(r)}$.
When summing over $z\in B(r)$, we bound according to the worse
scenario, choose $\delta$ small enough, and recall that
$\|l_m\|_{B(r)}\le m$,
\be{sub.10}
\begin{split}
\sum_{z\in B(r)} l_m^2(z)\pare{\frac{q_m(z)}{l_m(z)}}^2\ge&
(1-\delta) \|q_m\|^2_{B(r)}-2\delta\sum_{B(r)} l_m^2(z)\\
\ge& (1-\delta) \|q_m\|^2_{B(r)}-2\delta m^2\\
\ge& (\xi^2(1-\delta)-2\delta)m^2\ge \pare{(1-\epsilon)m}^2.
\end{split}
\ee

We can now state our subadditive result which we prove in the Appendix.
\bl{lem-subadd}
For $\xi>0$ small enough, for
any $r>0$ and for any $\delta$ small enough, there is
a constant $\J(\xi,r,\delta)$ such that
\be{sub-stat}
\lim_{m\to\infty} \frac{1}{m} \log\ P\pare{ \|\ind_{B(r)} l_m
\pi_{\delta}\pare{\frac{q_m(z)}{l_m(z)}}\|_2\ge\xi m}=\ -\J(\xi,r,\delta).
\ee
\el

\subsection{On the upper bound for the LDP.}\label{sec-UBLDP}
We show the following upper bound.
\bp{lem-LDUB} Assume $n^{2/3}\ll \xi_n\ll n^2$, and $d\ge 3$. Then,
\be{LDUB.1}
\begin{split}
\forall \epsilon>0,\quad \exists r_0>0,\quad \exists \alpha_0>0,\quad 
\exists \delta_0>0,&
\qquad \forall r>r_0,\quad \forall \alpha>\alpha_0, \quad \forall 
\delta<\delta_0\\
\lim\sup_{n\to\infty} \frac{1}{\sqrt{\xi_n}}\log\ P(X_n\ge \xi_n)\le &
-\alpha \J\pare{\frac{\sqrt{(1-\epsilon)}}{\alpha},r,\delta}+
\epsilon.
\end{split}
\ee
\ep
\bpr
The random walk cannot escape $[-n,n]^d$ in a time $n$.
Fix $\epsilon>0$. By Corollary~\ref{lem-D}, and 
at the expense of a polynomial term, there is a constant $M_A$,
and a finite volume $\Lambda_n$ with $|\Lambda_n|\le M_A$ such that
\be{sub.2}
P(X_n\ge \xi_n)\le n^\gamma \quad
P\pare{ \sum_{z\in \Lambda_n} q_n^2(z)\ge (1-\epsilon)\xi_n}.
\ee
Recalling that the walk is transient, it is convenient to
pass to an inifinite time-horizon. Thus, we define
\be{def-qinfinity}
l_{\infty}(z)=\sum_{i\in \N} \ind\acc{S(i)=z},\quad\text{and}\quad
q_{\infty}(z)=\sum_{i=1}^{l_{\infty}(z)} \eta_z(i).
\ee
We use now our monotony of the square charges, Corollary \ref{lem-mon},
to conclude
\be{sub.2infinity}
P(X_n\ge \xi_n)\le n^\gamma \quad
P\pare{ \sum_{z\in \Lambda_n} q^2_\infty(z)\ge (1-\epsilon)\xi_n}.
\ee
\cite{A07} establishes that $\Lambda_n$
can be transfered into a domain of finite diameter $\tilde \Lambda$.
Let $\T:\Lambda_n\to\tilde \Lambda$ be the {\it transfer map}
(see Proposition 5.2 of \cite{A07}). 
The following more precise statement is established in \cite{A07}: 
for any $\epsilon>0$, there is $r_0>0$, such that 
$\tilde \Lambda\subset B(r_0)$ for any large integer $n$,
and for any sequence of integers $\{k(z),z\in \Lambda_n\}$, with
$k(z)\le A\sqrt{\xi_n}$ for all $z\in \Lambda_n$, we have
\be{sub.3}
\begin{split}
P&\pare{ \sum_{z\in \Lambda_n} q_\infty^2(z)\ge (1-\epsilon)\xi_n\ 
l_\infty(z)=k(z),\ \forall z\in \Lambda_n}\\
\quad\quad& \le e^{\epsilon \sqrt{\xi_n}}\P_0\pare{
l_\infty(\T z)\ge k(z),\forall z\in \Lambda_n}Q\pare{ \sum_{z\in \Lambda_n}
\pare{\sum_{i=1}^{k(z)}\eta_z(i)}^2\ge (1-\epsilon)\xi_n}.
\end{split}
\ee
In the sum of the $\eta_z(i)$ over $[1,k(z)]$,
in the right hand side of \reff{sub.3},
we need to replace $k(z)$ by the larger value $l_\infty(\T z)$.
Again, we require a monotony of the $l_2$-norm of
the charges (i.e. Corollary~\ref{lem-mon} in
the Appendix), with the consequence that for a fixed realization
of the walk with $\{l_\infty(\T(z))\ge k(z),\forall z\in \Lambda_n\}$,
and any $r>r_0$, and with two shorthand notations
\[
A_n=(1-\epsilon)\xi_n,\quad\text{and}\quad
q_z(n)=\sum_{i=1}^n \eta_z(i).
\]
\be{sub.4}
\begin{split}
Q\pare{ \sum_{z\in \Lambda_n} q_z^2(k(z))\ge A_n}
\le& Q\pare{ \sum_{z\in \Lambda_n} q_z^2(l_\infty(\T z))\ge A_n}
= Q\pare{ \sum_{z\in \tilde \Lambda} q_z^2(l_\infty(z))\ge A_n}\\
\le& Q\pare{ \sum_{z\in B(r)}
q_z^2(l_\infty(z))\ge A_n}.
\end{split}
\ee
Thus, after averaging over the walk in \reff{sub.4}, and
summing over the $\{k(z),\ z\in \Lambda_n\}$
each term of \reff{sub.4}, we have that
for any $\epsilon>0$, there is $r>0$ such that
\be{sub.5}
\begin{split}
P&\pare{ \|\ind_{\Lambda_n} q_\infty\|^2_2\ge A_n}\le
e^{\epsilon \sqrt{\xi_n}}\sum_{{\bf k}}\E\cro{
\ind_{\acc{l_\infty(\T(z))\ge k(z),\forall z\in \Lambda_n}}
Q\pare{\sum_{z\in B(r)} q_z^2(l_\infty(z))\ge A_n}}\\
&\le e^{\epsilon \sqrt{\xi_n}}E\cro{ \prod_{z\in \tilde\Lambda} l_\infty(z)
Q\pare{\sum_{z\in B(r)} q_z^2(l_\infty(z))\ge A_n}}\\
&\le e^{\epsilon \sqrt{\xi_n}}(A\sqrt{\xi_n})^{|B(r)|}
P\pare{\sum_{z\in B(r)} q_z^2(l_\infty(z))\ge A_n,\
\max_{B(r)} l_\infty\le A\sqrt{\xi_n}}\\
&\quad
+e^{\epsilon \sqrt{\xi_n}}\E\cro{\prod_{z\in \tilde\Lambda} l_\infty(z)
\ind_{\{\exists z\in \tilde\Lambda,\ l_\infty(z)>A\sqrt{\xi_n}\}}}.
\end{split}
\ee
The second term in the right hand side of \reff{sub.5} is bounded
by a term $\exp(-\chi_d A\sqrt{\xi_n})$, whereas the first term
is estimated as follows.
\be{time-5}
\begin{split}
P\big(\sum_{z\in B(r)}q_z^2(l_\infty(z))&\ge A_n,\
\max_{B(r)} l_\infty\le A\sqrt{\xi_n}\big)\\
\le&(A\sqrt{\xi_n})^{|B(r)|} \sup_{k}
Q\pare{\sum_{z\in B(r)} q_z^2(k(z))\ge A_n}\
\P_0\pare{ l_{\infty}(z)=k(z),\ \forall z\in B(r)},
\end{split}
\ee
where the supremum in \reff{time-5} is over
integer sequences such that $\max_{B(r)} k(z)\le A\sqrt{\xi_n}$.

Now, we proceed similarly as in Section 6 of \cite{A07}.
We choose an integer sequence $\{ k(z),z\in B(r)\}$ with $
\max_{B(r)} k(z)\le A\sqrt{\xi_n}$, and define $|k|=\sum_{B(r)}k(z)$, and
\be{time-1}
\EE(k)=\acc{{\bf z}=(z(1),\dots,z(|k|))
\in \tilde\Lambda^{|k|}:\ \sum_{i=1}^{|k|} \ind\acc{z(i)=x} =k(x),\
\forall x\in B(r)}.
\ee
Now, for $\{k(z),z\in B(r)\}$ with $\max_{B(r)}k(z)\le A\sqrt{\xi_n}$,
we have, if $T_x=\inf\{n\ge 1: S_n=x\}$ and $T=\min\{T(x),x\in B(r)\}$
\be{time-2}
P_0\pare{l_{\infty}(z)=k(z),\ \forall z\in B(r)}=\sum_{{\bf z}\in \EE(k)}
\prod_{i=0}^{|k|-1}
P_{z(i)}\pare{T(z({i+1}))=T<\infty}P_{z(|k|)}(T=\infty).
\ee
We show now that if $d\ge 3$, and any $\epsilon>0$, there is
$\alpha(r,\epsilon)$ such that for all $x,y\in B(r)$
\be{time-3}
P_{x}\pare{T(y)=T<\infty}\le (1-\epsilon)
P_{x}\pare{T(y)=T<\alpha(r,\epsilon)}.
\ee
This would imply that
\be{time-4}
\begin{split}
P_0(l_{\infty}(z)=k(z),\ \forall z\in B(r))\le&
(1-\epsilon)^{|k|} \sum_{{\bf z}\in \EE(k)}
\prod_{i=0}^{|k|-1}
P_{z(i)}\pare{T(z({i+1}))=T<\alpha}P_{z(|k|)}(T=\infty)\\
\le & P_0(l_{\alpha(r,\epsilon) |k|}(z)\ge k(z),\ \forall z\in B(r)).
\end{split}
\ee
Using \reff{time-5}, \reff{time-4} and Corollary~\ref{lem-mon} in
the Appendix, we would obtain for $\alpha\ge |B(r)|\alpha(r,\epsilon)$
\be{time-7}
P\pare{\sum_{z\in B(r)}q_z^2(l_\infty(z))\ge A_n,\
\max_{B(r)} l_\infty\le A\sqrt{\xi_n}}\le (A\sqrt{\xi_n})^{|B(r)|}
P\pare{\sum_{z\in B(r)}q_z^2(l_{\alpha \sqrt{\xi_n}}(z))\ge A_n}.
\ee
We now prove \reff{time-3}. A classical result yields that for
aperiodic symmetric walk, and any positive integer $k$, we have $P_x(T_x=k)
\le c_d/k^{d/2}$ for some positive constant $c_d$.
Also, if we only consider pairs $x,y\in B(r)$ such
that $P_x(T_y=T<\infty)>0$, then there is an integer $l_r$ and
a positive constant $c_r$, such that
\be{time-8}
\inf_{x,y\in B(r)} P_x(T_y=T<l_r)=c_r.
\ee
Now, by conditioning
\be{time-9}
\frac{c_d}{(k+l_r)^{d/2}}\ge P_x(T_x=k+l_r)\ge
P_x(T=T_y=k)P_y(T_x=T=l_r).
\ee
Thus,
\be{time-10}
P_x(T=T_y=k)\le \frac{c_d}{c_r} \frac{1}{(k+l_r)^{d/2}}.
\ee
Therefore, there is $C>0$ such that
for any $x,y\in B(r)$ (with $P_x(T_y=T<\infty)>0$), and any integer $k$
\be{time-11}
P_x(k<T_y=T<\infty)=\sum_{i>k} P_x(T=T_y=i)\le \frac{C}{k^{d/2-1}}.
\ee
We conclude that for any $\epsilon$, there is $k(r,\epsilon)$ such
that for any $x,y\in B(r)$
\be{time-12}
P_x(k(r,\epsilon)<T_y=T<\infty)\le \epsilon.
\ee
\reff{time-12} implies that there is $\alpha(r,\epsilon)$ such that
\reff{time-3} holds.

The purpose of squeezing $\D_n^*(A)$ inside $B(r)$ is to renormalize time.
Indeed, the walk typically visits $\sqrt{\xi_n}$-times sites
of $\tilde\Lambda_n$ in a total time of order $\sqrt{\xi_n}$. 
Thus, section VII of
\cite{A07} establishes that there is $\alpha_0>0$ such that for
$\alpha>\alpha_0$, calling $m_n$ the integer part of $\alpha \sqrt{\xi_n}$,
we have for some $\gamma>0$
\be{sub.6}
P\pare{\|q_\infty\|_{B(r)}\ge \sqrt{(1-\epsilon)\xi_n}}\le
n^\gamma\ P\pare{\|q_{m_n}\|_{B(r)}
\ge m_n \frac{\sqrt{(1-\epsilon)}}{\alpha},S(m_n)=0}.
\ee
As in \cite{A07} (7.19), 
in order that the walk returns to the origin at time $m_n$,
we needed to add an piece of path
of arbitrary length $n$, satisfying $\{S(0)=S(n)=0\}$
whose probability is polynomial in $n$. Under our hypothesis of
aperiodicity of the walk, this latter fact is true.
Putting together \reff{sub.8}, \reff{sub.2}, \reff{sub.2infinity},
\reff{sub.6}, and invoking Lemma~\ref{lem-subadd}, we obtain
\reff{LDUB.1} and conclude our proof.
\epr

\subsection{Lower Bound in Theorem~\ref{theo-LDP}}\label{LB-LDP}
We first show a lower bound, similar to Proposition~\ref{lem-LDUB},
and then in Section~\ref{sec-limit}, we take a limit as
$\alpha,r$ go to infinity.
\bp{lem-LDLB} Assume $\frac{2}{3}<\beta<2$, and $d\ge 3$. Then,
\be{LDLB.1}
\forall \epsilon>0,\ \forall r,\ \forall \alpha, \ \forall \delta>0,\quad
\lim\inf_{n\to\infty} \frac{1}{\sqrt{\xi_n}}\log\ P(X_n\ge \xi_n)\ge 
-\alpha \J\pare{\frac{\sqrt{ \xi (1+\epsilon)}}{\alpha},r,\delta}.
\ee
\ep

In order to use Lemma~\ref{lem-subadd}, we need to
show that there is $C$ which might depend on $(r,\alpha,\delta)$ such that
(recall that $m_n$ is the integer part of $\alpha \sqrt{\xi_n}$)
\be{LB.aim}
P\pare{ \|q_n\|^2_{\Z^d} -n \ge \xi_n}\ge
C\ P\pare{\|\ind_{B(r)} l_{m_n}\pi_{\delta}(\frac{q_{m_n}}{l_{m_n}})
\|_2^2 \ge  (1+\epsilon) \xi_n,\ S(m_n)=0}.
\ee
Since $\pi_\delta(x)\le x$ for $x\ge 0$ and $\delta>0$, 
it is obvious that
\reff{LB.aim} follows from Lemma~\ref{lem-LB} with
the following choice. For a fixed $r>0$, we set $\D_n=B(r)$ and
$M_n=|B(r)|$. 
We take $m_n$ as the integer part of $\alpha \sqrt{\xi_n}$ for 
$\xi_n\ll n^2$,
with $\alpha$ as large as we wish, and $n$ going to infinity.
Note that from \reff{origin}, we have
\[
\P_0\pare{l_n(B(r))\ge \frac{\epsilon}{2} \xi_n}\le
\exp\pare{-\frac{\tilde \kappa_d \epsilon \xi_n}{2|B(r)|^{2/d}}},
\]
which is negligible. The term $\P_0(\|l_n\|_2^2\ge n^{2\beta-\epsilon})$
is dealt with Corollary~\ref{cor-2beta} since $\zeta_d(\beta)> \beta/2$.

\subsection{About the rate function}\label{sec-limit}
Using Lemmas \reff{lem-LDUB} and \reff{lem-LDLB}, we have
\be{synth.1}
\begin{split}
\forall \epsilon>0,\quad \exists r_0>0,\quad \exists \alpha_0>0,\quad
\exists \delta_0>0,&
\qquad \forall r,r'>r_0,\quad \forall \alpha,\alpha'>\alpha_0, 
\quad \forall \delta,\delta'<\delta_0\\
\alpha \J\pare{\frac{\sqrt{ \xi (1-\epsilon)}}{\alpha},r,\delta}
-\epsilon\le&
\alpha' \J\pare{\frac{\sqrt{ \xi (1+\epsilon)}}{\alpha'},r',\delta'}.
\end{split}
\ee
If we set 
\[
\v(x,r,\delta)=\frac{\J(x,r,\delta)}{x},
\]
we note that \reff{intro.5}, \reff{LDUB.1} and \reff{LDLB.1}
imply that $\v(x,r,\delta)$ is bounded as follows.
\be{synth.10}
\frac{c_+}{\sqrt{ 1+\epsilon}}\le \v(x,r,\delta)\le
\frac{1}{\sqrt{ 1-\epsilon}}\cro{c_-+\frac{\epsilon}{\sqrt \xi}}.
\ee
Now, \reff{synth.1} reads as
\be{synth.2}
\begin{split}
\forall \epsilon>0,\quad \exists r_0>0,\quad \exists x_0>0,\quad
\exists \delta_0>0,&
\qquad \forall r,r'>r_0,\quad \forall x,x'<x_0, 
\quad \forall \delta,\delta'<\delta_0\\
\v(x',r',\delta')\ge&{\sqrt{\frac{1-\epsilon}{1+\epsilon}}}
\v(x,r,\delta) -\frac{\epsilon}{\sqrt{ \xi (1+\epsilon)}}.
\end{split}
\ee
As we consider
subsequences when $x\to 0,\ r\to\infty,\ \delta\to 0$, we obtain
for any $\epsilon$ small
\be{synth.4}
\liminf \v\ge {\sqrt{\frac{1-\epsilon}{1+\epsilon}}}
\limsup\v -\frac{\epsilon}{\sqrt{\xi (1+\epsilon')}}
\ee
As $\epsilon$ vanishes in\reff{synth.4}, we conclude that
$\v(x,r,\delta)$ (along any subsequences) converges to a 
constant $\Q_2>0$.
\section{Appendix}\label{sec-appendix}
\subsection{Proof of Lemma~\ref{lem-russian}}
Note that since we assume that $\eta$ has a symmetric law,
\[
Q(\zeta(n)>t)=2Q\pare{\frac{1}{\sqrt n} \sum_{i=1}^n \eta(i)\ge \sqrt t}.
\]
We first treat the case $\eta\in \H_1$, and $\lambda_0>0$
is such that $E[\exp(\lambda_0 \eta)]<\infty$.
We use a Chebychev's exponential inequality. For $\lambda>0$,
\be{app.30}
\begin{split}
Q\pare{\frac{1}{\sqrt n} \sum_{i=1}^n \eta(i)\ge t}\le &
e^{-\lambda t} \pare{ E_Q\cro{ \exp\pare{\frac{\lambda}{\sqrt n}
\eta}}}^n\\
\le & e^{-\lambda t} \pare{ 1+\frac{\lambda^2}{2n}+
\frac{\lambda^3}{3! n \sqrt{n}}E\cro{|\eta|^3e^{\frac{\lambda}{\sqrt n}
|\eta|}}}^n.
\end{split}
\ee
First, choose $\lambda=\lambda_0 \sqrt{n}/2$. There is a constant 
$c_1$ such that
\be{app.31}
\begin{split}
Q\pare{\frac{1}{\sqrt n} \sum_{i=1}^n \eta(i)\ge t}\le &
\exp\pare{-\frac{\lambda_0\sqrt n}{2} t} 
\pare{ 1+\frac{\lambda_0^2}{8}+E\cro{\lambda_0^3|\eta|^3
e^{\lambda_0 |\eta|/2}}}^n\\
\le & \exp\pare{-\frac{\lambda_0\sqrt n}{2} t}
\pare{ 1+\frac{\lambda_0^2}{8}+c_1E\cro{ e^{\lambda_0 \eta}}}^n\\
\le & \exp\pare{-\frac{\lambda_0\sqrt n}{2} t+\beta_1 n},
\quad\text{with}\quad
\beta_1=\frac{\lambda_0^2}{8}+c_1E\cro{ e^{\lambda_0 \eta}}.
\end{split}
\ee
Let $\beta_0=4\beta_1/\lambda_0$ and note that for 
$t\ge \sqrt{\beta_0 n}$, we have
\be{app.32}
Q\pare{\frac{1}{\sqrt n} \sum_{i=1}^n \eta(i)\ge t}\le 
\exp(-\frac{\lambda_0}{4} \sqrt{n} t).
\ee
Now, we assume that $t\le \sqrt{\beta_0 n}$, and we choose
$\lambda=\gamma t$ for $\gamma$ to be adjusted latter. Inequality
\reff{app.30} yields
\be{app.33}
Q\pare{\frac{1}{\sqrt n} \sum_{i=1}^n \eta(i)\ge t}\le 
e^{-\gamma t^2}\pare{ 1+\frac{\gamma^2 t^2}{2n}+
\frac{\gamma^3 t^3}{n\sqrt n}E\cro{ |\eta|^3
\exp\pare{ \frac{\gamma t}{\sqrt n}|\eta|}}}^n.
\ee
$\gamma$ needs to satisfy many constraints. First, in order for the
exponential of $|\eta|$ in \reff{app.33} to be finite, we need
that $\gamma \sqrt{\beta_0}\le \lambda_0/2$, in which case
\be{app.34}
\begin{split}
Q\pare{\frac{1}{\sqrt n} \sum_{i=1}^n \eta(i)\ge t}\le &
e^{-\gamma t^2}\pare{ 1+\frac{\gamma^2 t^2}{2n}+
\frac{\gamma^3 t^3}{\lambda_0^3 n\sqrt n} c_1 E\cro{ e^{\lambda_0\eta}}}^n\\
\le &\exp\pare{-\gamma t^2+\frac{\gamma^2 t^2}{2}+
\frac{\gamma^3 t^3}{\lambda_0^3 \sqrt n} c_1 E\cro{ e^{\lambda_0\eta}}}.
\end{split}
\ee
In the right hand side of \reff{app.34}, the term in $\gamma^2$ is
innocuous as soon as $\gamma\le 1/2$. Also, since $t\le \sqrt{\beta_0 n}$,
the term in $\gamma^3$ is innocuous as soon as 
\[
\gamma\le \frac{\lambda_0^3}{\sqrt{\beta_0}c_1 E\cro{ e^{\lambda_0\eta}}},
\]
so that \reff{app.34} yields
\be{app.35}
Q\pare{\frac{1}{\sqrt n} \sum_{i=1}^n \eta(i)\ge t}\le
\exp(-\kappa_\infty t^2),\quad\text{with}\quad
\kappa_\infty=\min\pare{\frac{\lambda_0^3}{\sqrt{\beta_0}c_1
E\cro{ e^{\lambda_0\eta}}},\frac{1}{2},\frac{
\lambda_0}{2\sqrt\beta_0}}.
\ee
Using \reff{app.32} and \reff{app.35}, it is immediate to
deduce \reff{prel.1}.

Assume now $\eta\in \H_\alpha$ for $1<\alpha<2$. From
Kasahara's Tauberian theorem, 
there is a constant $\kappa_\alpha$
and $\beta_0>0$ such 
that for $t\ge \sqrt{\beta_0 n}$, we have 
\be{app.36}
Q\pare{\frac{1}{\sqrt n} \sum_{i=1}^n \eta(i)\ge t}\le 
\exp\pare{-\kappa_\alpha \frac{(t\sqrt n)^\alpha}{n^{\alpha-1}}}.
\ee
Now, when $t\le \sqrt{\beta_0 n}$, we use the argument of
the previous case $\H_1$, to obtain \reff{app.35}.

Finally, when $\eta\in \H_2$, Chen has shown in \cite{chen07} 
there is a constant $C$ such that for any $k\in \N$
\be{chen1}
E_Q\cro{\pare{ \sum_{i=1}^n \eta(i)}^{2k}}\le C^k k! n^k.
\ee
\reff{chen1} implies that for some $\lambda_1>0$
\be{chen2}
\sup_n E_Q\cro{\exp\pare{ \lambda_1\pare{ \frac{\sum_{i=1}^n \eta(i)}
{\sqrt n}}^2}}<\infty.
\ee
Thus, there is a constant $C_1$ such that for any $n\in \N$, and $t>0$
\be{chen3}
P_Q(\zeta(n)>t)\le C_1 e^{-\lambda_1 t}.
\ee
This concludes the proof of Lemma~\ref{lem-russian}.

\subsection{Proof of Lemma~\ref{lem-subadd}}
We fix two integers $k$ and $n$, with $k$ to be taken first
to infinity. Let $m,s$ be integers such that $k=mn+s$, and $0\le s<n$.
The phenomenon behind the subadditive argument (to come) is that
the rare event
\[
\A_k(\xi,r,\delta)=\acc{\|l_k \pi_{\delta}
(\frac{q_k}{l_k})\|_{B(r)}\ge \xi k,\quad S(k)=0},
\]
can be built by concatenating the ${\it same }$ optimal scenario
realizing $\A_n(\xi,r,\delta)$ on $m$ consecutive periods of length
$n$, and one last period of length $s$ where the scenario is
necessarly special and its cost innocuous. The crucial independence
between the different period is obtained by forcing the walk
to return to the origin at the end of each period.

Thus, our first step is to exhibit an optimal strategy
realizing $\A_n(\xi,r,\delta)$. For this purpose,
we show that a finite number of values of the discrete variables
$\{l_n(z), \pi_{\delta}(q_n/l_n(z)),\ z\in B(r)\}$ are needed
to estimate the probability of $\A_n(\xi,r,\delta)$.

First recall that $\xi$ is as small as we wish. In particular, we take
it such that $Q(\eta>\xi)>0$.
Note that $\sum_{B(r)}l_n(z)\le n$. 
Also, using Lemma~\ref{lem-russian}, there is a constant $C$ such that 
\[
P\pare{\exists z\in B(r),\ \zeta_z(l_n(z))\ge An}\le \exp(-CAn).
\]
Consequently, the same holds for
$q_n(z)/l_n(z)=\sqrt{\zeta_z(l_n(z))/l_n(z)}$. Thus, if we denote by
\be{sub.20}
\EE(A)=\acc{\exists z\in B(r),\ \frac{q_n(z)}{l_n(z)}\ge \sqrt{An}},
\quad\text{then}\quad P(\EE(A))\le \exp(-CAn).
\ee

On the other hand, there is an obvious lower bound obtained by considering 
monomers making up one single pile:
\be{sub.24}
\S=\acc{l_n(0)=n,\ \eta(i)>\xi,\ \forall i<n\ },\quad\text{and}
\quad P(\S)=\pare{\frac{1}{2d+1}}^n\ Q(\eta>\xi)^n.
\ee
Thus, using \reff{sub.20} and \reff{sub.24}, we have for $A$ large enough
\be{sub.21}
2 P(\EE(A))\le P\pare{ \A_n(\xi,r,\delta)}
\Longrightarrow P\pare{ \A_n(\xi,r,\delta)} \le
2P\pare{ \A_n(\xi,r,\delta) ,\ \EE^c(A)}.
\ee
We conclude that for $z\in B(r)$ there are
\be{sub.23}
\lambda_n(z)\in [0,n]\cap \N,\quad\text{and}\quad 
\kappa_n(z)\in [0,\sqrt{An}]\cap \delta \N,\quad\text{with}\quad
\|\lambda_n \kappa_n\|_{B(r)}\ge \xi n,
\ee
such that
\be{sub.22}
P\pare{ \A_n(\xi,r,\delta)}\le \pare{\frac{2n\sqrt{An}}{\delta}}^{|B(r)|}
P\pare{ l_n|_{B(r)}=\lambda_n,\ 
\pi_{\delta}(\frac{q_n}{l_n})|_{B(r)}=\kappa_n}.
\ee
The factor 2 in the constant appearing in the right hand side
of \reff{sub.22} is to account for the choice of a positive total
charge on all sites of $B(r)$, by using the symmetry of the charge's
distribution.
Let $z^*\in B(r)$ be a site where 
\be{sub.25}
\lambda_n(z^*) \kappa^2_n(z^*)=\max_{B(r)} \lambda_n \kappa^2_n,
\quad\text{and note that}\quad
\lambda_n(z^*) \kappa^2_n(z^*)\ge n\xi^2.
\ee
Indeed, one uses that $\sum_{B(r)}\lambda_n(z)\le n$ and
\[
(n\xi)^2\le \sum_{z\in B(r)} (\lambda_n(z)\kappa_n(z))^2\le
\max_{B(r)}\pare{\lambda_n \kappa^2_n}\sum_{z\in B(r)} \lambda_n(z)
\le n\lambda_n(z^*) \kappa^2_n(z^*).
\]
Also,
\be{sub.26}
P\pare{\pi_\delta\pare{ \frac{q_n(z^*)}{l_n(z^*)}}=\kappa_n(z^*)}>0
\Longleftarrow
P\pare{\eta\ge  \kappa_n(z^*)}>0.
\ee
We set 
$\lambda_r(z^*)=r$, and $\lambda_r(z)=0$ for $z\not= z^*$. We define
the following symbols, for integers $i<j$
\be{sub.27}
l_{[i,j[}(z)=\sum_{t=i}^{j-1} \ind\acc{S_t-S_{i-1}=z},\quad\text{and}
\quad
q_{[i,j[}(z)=\sum_{t=i}^{j-1} \eta(t)\ind\acc{S_t-S_{i-1}=z}.
\ee
Note that on disjoints sets $I_k=[i_k,j_k[$, the variables
$\{(l_{I_k},q_{I_k}),k\in \N\}$ are independent. Finally, we define
and the following sets, for $i=1,\dots,m$
\be{sub.28}
\begin{split}
\A_n^{(i)}=&\acc{ l_{[(i-1)n,in[}|_{B(r)}=\lambda_n,\ 
\pi_{\delta}\pare{\frac{q_{[(i-1)n,in[}}{l_{[(i-1)n,in[}}}|_{B(r)}=
\kappa_n,\ S_{in}=0}\\
\A_s=&\acc{l_{[mn,k[}=s \delta_{z^*},\ \eta(i)>\kappa_n(z^*),
\forall i\in [mn,k[}.
\end{split}
\ee

On the event $\{S_{in}=0,\ \forall i=0,\dots,m\}$, the local charges,
and local times in $[0,k[$ on site $z$ are respectively
\be{sub.29}
q_k(z)=\sum_{i=1}^m q_{[(i-1)n,in[}(z)+q_{[mn,k[}(z),\quad\text{and}\quad
l_k(z)=\sum_{i=1}^m l_{[(i-1)n,in[}(z)+l_{[mn,k[}(z).
\ee
We need to show now that
\be{sub.30}
\bigcap_{i=1}^m \A_n^{(i)}\cap\A_s\subset \A_k(\xi,r,\delta).
\ee
In other words, we need to see that under $\bigcap_i\A_n^{(i)}\cap\A_s$,
we have
\be{sub.31}
\sum_{z\in B(r)} \pare{m\lambda_n(z)+\lambda_s(z)}^2
\pi^2_{\delta}\pare{ \frac{q_k(z)}{m\lambda_n(z)+\lambda_s(z)}}
\ge (\xi k)^2.
\ee
Note that under $\cap_i\A_n^{(i)}\cap\A_s$, for any $z\in B(r)$
\be{sub.32}
\frac{q_{[(i-1)n,in[}(z)}{l_{[(i-1)n,in[}(z)}\in [\kappa_n(z),
\kappa_n(z)+\delta[\Longrightarrow
\frac{\sum_{i=1}^m q_{[(i-1)n,in[}(z)}{ml_{[(i-1)n,in[}(z)}
\in [\kappa_n(z),\kappa_n(z)+\delta[.
\ee
Thus, if $s=0$, \reff{sub.31} would hold trivially.

We assume for simplicity that $z^*=0$, and postpone to Remark~\ref{rem-z}
the general case. Note that for $z=z^*=0$
\be{sub.33}
\frac{q_{[mn,k[(0)}}{\lambda_s(0)}\ge \kappa_n(0)\quad\text{and 
\reff{sub.32} imply that}\quad
\frac{\sum_{i=1}^m q_{[(i-1)n,in[}(0)+q_{[mn,k[(0)}}
{m\lambda_n(0)+s}\ge  \kappa_n(0),
\ee
whereas for $z\not=0$, $q_k(z)=q_{mn}(z)$ and $l_k(z)=l_{mn}(z)$, so
that checking \reff{sub.31} reduces to checking
\be{sub.34}
\begin{split}
\pare{ m \lambda_n(0)+s}^2& \pi^2_\delta\pare{
\frac{\sum_{i=1}^m q_{[(i-1)n,in[}(0)+q_{[mn,k[(0)}}
{m\lambda_n(0)+s}}\\
&\qquad\qquad 
-\pare{ m \lambda_n(0)}^2 \pi^2_\delta\pare{
\frac{\sum_{i=1}^m q_{[(i-1)n,in[}(0)}{m\lambda_n(0)}}
\le\pare{k^2-(mn)^2}\xi^2.
\end{split}
\ee
Using \reff{sub.33}, it is enough to check that
\be{sub.35}
\pare{ 2m\lambda_n(0)+s} \kappa_n(0)^2\ge
(2mn+s) \xi^2.
\ee
Recall that \reff{sub.25} yields
$\lambda_n(z^*) \kappa^2_n(z^*)\ge n\xi^2$ so that when $z^*=0$
\be{sub.37}
\lambda_n(0)\kappa_n(0)^2\ge n \xi^2,\quad\text{and}\quad
\kappa_n(0)^2\ge \xi^2,
\ee
which implies \reff{sub.35} right away.

Now, \reff{sub.30} implies that for $c,c'$ depending on $\delta,r$ and
$A$, we have
\be{sub.38}
\begin{split}
P\pare{ \A_n(\xi,r,\delta)}^m P(\A_s)&\le
(cn)^{c'm} P(\A_n^{(1)})\dots P(\A_n^{(m)}) P(\A_s)\\
&\le (cn)^{c'm} P\pare{ \bigcap_{i\le m}\A_n^{(i)}\cap \A_s}\\
&\le (cn)^{c'm} P\pare{ \A_k(\xi,r,\delta)}.
\end{split}
\ee

We now take the logarithm on each side of \reff{sub.38}
\be{sub.39}
\frac{nm}{nm+s} \frac{\log(P(\A_n(\xi,r,\delta)))}{n}+  
\frac{\log(P(\A_s))}{k}\le
\frac{c'm\log(c n)}{nm+s}+\frac{\log(P(\A_k(\xi,r,\delta)))}{k}.
\ee
We take now the limit $k\to\infty$ while $n$ is kept fixed 
(e.g. $m\to\infty$) so that
\be{sub.40}
\frac{\log(P(\A_n(\xi,r,\delta)))}{n}\le \frac{c'\log(c n)}{n}
+\liminf_{k\to\infty} \frac{\log(P(\A_k(\xi,r,\delta)))}{k}.
\ee
By taking the limit sup in \reff{sub.40} as $n\to\infty$, 
we conclude that the limit in \reff{sub-stat} exists.

\br{rem-z} We treat here the case $z^*\not= 0$. Note that
this is related to the strategy on a single period of length $s$.
If we could have that monomers in a piece of length $s$
pile up on site $z^*$, then \reff{sub.30} would hold since it only
uses that $\lambda_n(z^*) \kappa^2_n(z^*)\ge n\xi^2$.
However, the walk starts at the origin, and
each period of length $n$ sees the walk returning to the origin.
The idea is to insert a period of length $s$ into
the first time-period of length $n$ at the first time the walk hits
$z^*$. Then, since we still use
a scenario with a single pile at site $z^*$ with charges exceeding $\xi$,
we should have
\be{sub.41}
l_{n+s}=\lambda_n+s\delta_{z^*},\ 
\pi_\delta(\frac{q_{n+s}}{l_{n+s}})|_{B(r)}
\ge  \kappa_n,\quad\text{and}\quad S_{n+s}=0.
\ee
More precisely, let $\tau^*=\inf\{n\ge 0:\ S_n=z^*\}$, and note that 
\be{sub.42}
P\pare{\A_n^{(1)}}=\sum_{i=0}^{n-1} P\pare{ \tau^*=i,\ \A_n^{(1)}}.
\ee
Let $i^*<n$ be such that 
\[
P\pare{ \tau^*=i^*,\ \A_n^{(1)}}=\max_{i<n} P\pare{ \tau^*=i,\ \A_n^{(1)}}.
\]
Then,
\be{sub.43}
P\pare{\A_n^{(1)}}\le n P\pare{ \tau^*=i^*,\ \A_n^{(1)}},
\ee
and, adding a subscript to $P$ to explicit the starting point of the
walk
\be{sub.44}
\begin{split}
P_0\pare{\A_n^{(1)}} &
P_{0}\pare{l_{[0,s[}=s \delta_{0},\ \eta(i)>\kappa_n(z^*),
\forall i\in [0,s[}\\
&\le n\ P_0\pare{ \tau^*=i^*,\ \A_n^{(1)}}
P_{z^*}\pare{l_{[0,s[}=s \delta_{z^*},\ \eta(i)>\kappa_n(z^*)}\\
&\le P_0\pare{ l_{[0,s[}=s ,\ \pi_{\delta}
\pare{\frac{q_{n+s}}{l_{n+s}}}|_{B(r)} \ge  \kappa_n,\ S_{n+s}=0}.
\end{split}
\ee
\er

\subsection{On a monotony property}
We prove in this section the following result which is a corollary of
Lemma 5.3 of \cite{AC04}. 
\bc{lem-mon} For any integer $n$, assume that $\{\eta_j(i),j=1,\dots,n,\ 
i\in \N\}$ are independent symmetric variables, and for any sequence
$\{n_j,n_j',j=1,\dots,n\}$ with $n_j'\ge n_j$, and any $\xi>0$, we have
\be{mon.1}
P\pare{ \sum_{j=1}^n\pare{\sum_{i=1}^{n_j} \eta_j(i)}^2>\xi}\le
P\pare{ \sum_{j=1}^n\pare{\sum_{i=1}^{n_j'} \eta_j(i)}^2>\xi}.
\ee
\ec
\bpr
We prove the result by induction. First, for $n=1$, we use
first the symmetry of the distribution of the $\eta$'s and
then Lemma 5.3 of \cite{AC04} to have for $n_1'\ge n_1$
\be{mon.2}
\begin{split}
P\pare{\pare{ \sum_{i=1}^{n_1} \eta_1(i)}^2>\xi}=&2
P\pare{\sum_{i=1}^{n_1} \eta_1(i)>\sqrt{\xi}}\le 2
P\pare{\sum_{i=1}^{n_1'} \eta_1(i)>\sqrt{\xi}}\\
\le& P\pare{\pare{ \sum_{i=1}^{n_1'} \eta_1(i)}^2>\xi}.
\end{split}
\ee
Now, assume that \reff{mon.1} is true for $n-1$, and
call $\Gamma_j=(\eta_j(1)+\dots+\eta_j(n_j))^2$ and $\Gamma_j'$ the 
sum of the $\eta_j$ up to
$n_j'$. We only write the proof in the case where $\Gamma_j$ has a density, 
say $g_{\Gamma_j}$. The case of a discrete distribution is trivially
adapted. Then
\be{mon.3}
\begin{split}
P\pare{\sum_{j=1}^n \Gamma_j>\xi}=& P(\Gamma_1>\xi)+
\int_0^{\xi} g_{\Gamma_1}(z) P\pare{\sum_{j=2}^n \Gamma_j>\xi-z}dz\\
\le & P(\Gamma_1>\xi)+ 
\int_0^{\xi} g_{\Gamma_1}(z) P\pare{\sum_{j=2}^n \Gamma_j'>\xi-z}dz\\
=&P\pare{\Gamma_1+\sum_{j=2}^n \Gamma_j'>\xi}
\end{split}
\ee
Then, we rewrite the sum on the right hand side of \reff{mon.3}
$\Gamma_1+(\Gamma_2'+\dots+\Gamma_n')=\Gamma_2'+(\Gamma_1+\dots+\Gamma_n')$, and single out $\Gamma_2'$
in the first step of \reff{mon.3} to conclude.
\epr

\end{document}